# The fractional orthogonal derivative.

Enno Diekema
e.diekema@gmail.com

**Abstract**.

This paper builds on the notion of the so-called orthogonal derivative, where an n-th order derivative is approximated by an integral involving an orthogonal polynomial of degree n. This notion was reviewed in great detail in a paper in J. Approx. Theory (2012) by the author and Koornwinder. Here an approximation of the Weyl or Riemann-Liouville fractional derivative is considered by replacing the n-th derivative by its approximation in the formula for the fractional derivative. In the case of, for instance, Jacobi polynomials an explicit formula for the kernel of this approximate fractional derivative can be given. Next we consider the fractional derivative as a filter and compute the transfer function in the continuous case for the Jacobi polynomials and in the discrete case for the Hahn polynomials. The transfer function in the Jacobi case is a confluent hypergeometric function. A different approach is discussed which starts with this explicit transfer function and then obtains the approximate fractional derivative by taking the inverse Fourier transform. The theory is finally illustrated with an application of a fractional differentiating filter. In particular, graphs are presented of the absolute value of the modulus of the transfer function. These make clear that for a good insight in the behavior of a fractional differentiating filter one has to look for the modulus of its transfer function in a log-log plot, rather than for plots in the time domain.

## 1    Introduction.

As soon as the differential calculus was invented by Newton and Leibniz, the question arose about the meaning of a derivative of non-integer order. Many mathematicians have worked at the subject and many different definitions were proposed.

An excellent historical survey is given in the book by Miller and Ross (1993). Also the papers of Osler (1970a,b), (1971), (1972a,b,c,d), (1973) and (1974) give many historical details. In the encyclopedic work of Samko et al. (1993) most of the usual definitions are given. The books of Kilbas et al. (2006) gives a survey of a number of definitions of the fractional calculus. The book of Ortigueira (2011) considers only fractional derivatives of causal functions. Last but not least, Tenreiro Machado et al. (2011) give a survey of the developments in the field from 1974.

Nowadays fractional derivatives are much used in fractional differential equations, which are very suitable for describing certain physical phenomena, better than when using ordinary derivatives. For instance, a fractional PID (Proportional, Integrating and Differentiating) controller can be made much simpler and more robust than an ordinary PID controller. Some other areas of application of the fractional derivatives are viscoelastic materials, hydrodynamics, rheology, diffusive transport, electrical networks, control theory, electromagnetic theory, signal and image processing and probability. The definition of fractional derivative may vary with the application, since in each case one will look for the most suitable definition.

Unlike with the ordinary derivative, there is no simple geometric description of the fractional derivative. However, see Podlubny (2002) and Tenreiro Machado (2009b).

The most common definitions of fractional derivative, which we also use in this paper, are by taking ordinary derivatives of a Riemann-Liouville or Weyl fractional integral. Even here, different choices can be made which give quite different results. For instance, with $v > 0$ and $\mu = n - v > 0$ we can write

(1.1) $\qquad \dfrac{d^v}{dx^v} e^x = \dfrac{d^n}{dx^n} \dfrac{1}{\Gamma(\mu)} \int_{-\infty}^{x} (x-y)^{\mu-1} e^y dy = e^x = \sum_{k=0}^{\infty} \dfrac{x^k}{\Gamma(k+1)}$

and thus obtain a natural choice for the fractional derivative of $e^x$, but we can also choose the definition

1.2|) $\qquad \dfrac{d^v}{dx^v} e^x = \dfrac{d^n}{dx^n} \dfrac{1}{\Gamma(\mu)} \int_{0}^{x} (x-y)^{\mu-1} e^y dy = \sum_{k=0}^{\infty} \dfrac{x^{k-v}}{\Gamma(k-v+1)}$

Here the right side follows from

(1.3) $\qquad \dfrac{d^v}{dx^v} x^k = \dfrac{d^n}{dx^n} \dfrac{1}{\Gamma(\mu)} \int_{0}^{x} (x-y)^{\mu-1} y^k dy = \dfrac{\Gamma(k+1)}{\Gamma(k-v+1)} x^{k-v}$

which looks like a natural choice for the fractional derivative of $x^k$. (Note that replacement of the lower bound of the last integral by $-\infty$ would give a divergent integral). Thus one has to be careful with just extrapolating a naturally looking definition of fractional derivative of some elementary functions.

In Diekema and Koornwinder (2012) we reviewed a formula with a long history for the so-called orthogonal derivative. This derivative (of order $n$) can be computed as a limit of a certain integral. This is a generalization of the usual notion of the $n$-th order derivative and, more important, when ignoring the limit it can be used as an approximation of the $n$-th order derivative. In the present paper the orthogonal derivative will be used in order to generalize the definitions of Riemann-Liouville and Weyl for the fractional derivative and to approximate these fractional derivatives. In the case of the orthogonal derivative associated with the Jacobi polynomials, the kernel of the resulting integral transform approximating the fractional derivative can be computed explicitly. In the case of the Hahn polynomials similar explicit results are obtained for the approximation of the fractional difference.

In general, Fourier and Laplace transforms are important tools for finding solutions of fractional orthogonal differential equations. So in our paper we also consider the action of the various operators in the frequency domain. We then consider the fractional derivative as a filter with a so-called transfer function. Instead of a convolution we then have multiplication of transfer functions. A picture of the modulus of the transfer function gives a very good insight how well the approximating operators behave.

The idea of combining an orthogonal derivative with a fractional integral in order to approximate fractional derivatives is already discussed in Chen and Chen (2011) and Liu et al (2012). However explicit formula like ours are not given in these papers. Also a discussion of the filter in the frequency domain is missing.

Let us now give a summary of the sections of this paper.

In section 2 the definitions and basic properties of the Riemann-Liouville and Weyl fractional integral and derivative are recalled. Also their Fourier transform are mentioned.

In section 3 the approximate Weyl and Riemann-Liouville fractional derivatives are defined by using the approximate $n$-th order derivative coming from the orthogonal derivative.

In section 4 we give explicit expressions in the case of the Jacobi polynomials. These simplify in the case of Gegenbauer and Legendre polynomials.

In section 5 we consider the transfer function of the approximate fractional derivative. Explicit results in the continuous case are given for Jacobi polynomials and in the discrete case for Hahn polynomials.

The general case suggests an alternative approach to obtain an approximate fractional derivative, by starting with a transfer function which approximates the transfer function of the fractional derivative and then taking the inverse Fourier transform.

This method is discussed in Section 6 for a transfer function given by a confluent hypergeometric function and in Section 7 for two more examples involving elementary transfer functions.

In section 8 the theory of the fractional derivative will be applied with a discrete function and so the fractional difference is used as a base for deriving a formula for the fractional orthogonal difference.

In section 9 the theory of the fractional difference is applied with the discrete Hahn polynomials.

In section 10 we derive the transfer function for the fractional difference with the discrete Hahn polynomials.

In section 11 the theory will be applied to a fractional differentiating filter. Analog as well discrete filters are treated. The theory is illustrated by the graphs of the modulus of the transfer functions.

## 2 The Riemann-Liouville and the Weyl fractional transforms.

We define two versions of the fractional integral: the Riemann-Liouville and the Weyl transform. In terms of these we can define a fractional derivative.

Let $\mu \in \mathbb{C}$ with $\operatorname{Re} \mu > 0$. Let $f$ be a function on $(-\infty, b)$ which is integrable on bounded subintervals. Then the Riemann-Liouville integral of order $\mu$ is defined as

(2.1) $\quad R^{-\mu}[f](x) = f^{(-\mu)}(x) := \dfrac{1}{\Gamma(\mu)} \displaystyle\int_{-\infty}^{x} f(y)(x-y)^{\mu-1} dy$

A sufficient condition for the absolute convergence of this integral is $f(-x) = O(x^{-\mu-\varepsilon}), \varepsilon > 0, x \to \infty$.

Similarly, for $\operatorname{Re} \mu > 0$ and $f$ a locally integrable function on $(a, \infty)$ such that $f(x) = O(x^{-\mu-\varepsilon})$, $\varepsilon > 0$, $x \to \infty$, the Weyl integral of order $\mu$ is defined as

(2.2) $\quad W^{-\mu}[f](x) = f^{(-\mu)}(x) := \dfrac{1}{\Gamma(\mu)} \displaystyle\int_{x}^{\infty} f(y)(y-x)^{\mu-1} dy$

Clearly

(2.3) $\quad R^{-\mu}[f(-\,\cdot\,)](-x) = W^{-\mu}[f](x)$

The Riemann-Liouville integral is often given for a function $f$ on $[0, b)$ in the form

$$(2.4) \quad R^{-\mu}[f](x) = f^{(-\mu)}(x) := \frac{1}{\Gamma(\mu)} \int_0^x f(y)(x-y)^{\mu-1} dy \quad (x > 0)$$

This can be obtained from (2.1) if we assume that $f(x) = 0$ for $x \leq 0$ (then $f$ is called a causal function). For instance, if $f(x) = x^\alpha$ $(x > 0)$ and $f(x) = 0$ $(x \leq 0)$ then

$$(2.5a) \quad R^{-\mu}[f](x) = \frac{\Gamma(\alpha+1)}{\Gamma(\alpha+\mu+1)} x^{\alpha+\mu} \quad (x > 0, \ \text{Re}(\alpha) > -1)$$

and

$$(2.5b) \quad W^{-\mu}[f](x) = \frac{\Gamma(-\alpha-\mu)}{\Gamma(-\alpha)} x^{\alpha+\mu} \quad (x > 0, \ \text{Re}(\alpha+\mu) < 0)$$

When $-1 < \text{Re}\,\alpha < \text{Re}(-\mu) < 0$ both formulas are valid. The quotients of Gamma functions are related by

$$(2.6) \quad \frac{\Gamma(-\alpha-\mu)}{\Gamma(-\alpha)} = \frac{\sin(\alpha\pi)}{\sin(\alpha+\mu)\pi} \frac{\Gamma(\alpha+1)}{\Gamma(\alpha+\mu+1)}$$

Substitution in (2.5) gives

$$(2.7) \quad W^{-\mu}[f](x) = \frac{\sin(\alpha\pi)}{\sin(\alpha+\mu)\pi} R^{-\mu}[f](x)$$

The fractional derivatives are defined by the following formulas with $\text{Re}\,\nu < n$, $n$ a positive integer and $D = \frac{d}{dx}$

$$(2.8a) \quad R^\nu[f](x) = D^n\big[R^{\nu-n}[f]\big](x)$$

$$(2.8b) \quad W^\nu[f](x) = (-1)^n D^n\big[W^{\nu-n}[f]\big](x)$$

Put $W^0 = id = R^0$. Under assumption of sufficient differentiability and convergence we have $W^\mu W^\nu = W^{\mu+\nu}$, $R^\mu R^\nu = R^{\mu+\nu}$ for all $\mu, \nu \in \mathbb{C}$ and $R^n = D^n$, $W^n = (-1)^n D^n$.

For $\text{Re}\,\nu < 0$ formula (2.8a) and (2.8b) are easily derived and for $\text{Re}\,\nu > 0$ they are definitions. This is under the condition that all derivatives of order less then $n+1$ of the function $f(x)$ should exist.

Because we need the Fourier transform of formulas (2.8) we use the following theorem proved by Samko et al. (1993) Ch.2 §7 for the Riemann-Liouville derivative as well for the Weyl derivative.

Theorem 2.1. Let $f$ and $g$ be functions on $\mathbb{R}$ for which the Fourier transforms exist and are given by

$$(2.9a) \quad F(\omega) = \frac{1}{\sqrt{2\pi}} \int_{-\infty}^\infty e^{i\omega x} f(x) dx$$

$$(2.9b) \quad G(\omega) = \frac{1}{\sqrt{2\pi}} \int_{-\infty}^\infty e^{i\omega x} g(x) dx$$

If $g = R^\nu[f]$ then the following relation holds:

$$(2.10) \quad \frac{G(\omega)}{F(\omega)} = (-i\omega)^\nu$$

If $g = W^\nu[f]$ then

$$(2.11) \quad \frac{G(\omega)}{F(\omega)} = (i\omega)^\nu$$

Here and in the rest of this paper $z^\nu = e^{\nu \ln z}$ with $-\pi < \arg(z) < \pi$. So $(-i\omega)^\nu = (-i\omega + 0)^\nu = e^{-\pi i \nu/2}(\omega + i0)^\nu$. If $z$ is on the cut $(-\infty, 0]$ and $z \neq 0$ then we distinguish between $(z+i0)^\nu = e^{i\pi\nu}(-z)^\nu$ and $(z-i0)^\nu = e^{-i\pi\nu}(-z)^\nu$.

The quotient $H(\omega) = \frac{G(\omega)}{F(\omega)}$ will be called the *transfer function*, where we follow the terminology of filter theory.

## 3  The fractional Weyl transform for the orthogonal derivative.

In Diekema and Koornwinder (2012) the following theorem is proved:

Theorem 3.1. Let $n$ be a positive integer. Let $p_n$ be an orthogonal polynomial of degree $n$ with respect to the orthogonality measure $\mu$ for which all moments exist. Let $x \in \mathbb{R}$. Let $I$ be a closed interval such that, for some $\varepsilon > 0$, $x + \delta v \in I$ if $0 \leq \delta \leq \varepsilon$ and $v \in \mathrm{supp}(\mu)$. Let $f$ be a continuous function on $I$ such that its derivatives of order $1, 2, \ldots, n$ at $x$ exist. In addition, if $I$ is unbounded, assume that $f$ is of at most polynomial growth on $I$. Then

(3.1a) $\qquad f^{(n)}(x) = \lim_{\delta \downarrow 0} D_\delta^n[f](x)$

where

(3.1b) $\qquad D_\delta^n[f](x) = \dfrac{k_n n!}{h_n} \dfrac{1}{\delta^n} \displaystyle\int_{\mathbb{R}} f(x + \delta \xi) p_n(\xi) d\mu(\xi)$

By the assumptions the integral converges absolutely. $h_n$ and $k_n$ are defined by

(3.2a) $\qquad h_n = \displaystyle\int_{\mathbb{R}} p_n(\xi)^2 d\mu(\xi)$

and

(3.2b) $\qquad p_n(x) = k_n x^n + q_{n-1}(x) \qquad$ with $q_{n-1}$ a polynomial of degree $< n$.

We call the limit given by the right-hand side of (3.1a) the *orthogonal derivative* at $x$ of order $n$. See Remark 3.3 for some weaker conditions.

Theorem 3.1. suggests the definition of the approximate fractional orthogonal derivative. We can apply formula (3.1a) in formula (2.8b). This yields

(3.3) $\qquad W^\nu[f](x) = \dfrac{(-1)^n}{\Gamma(n-\nu)} \dfrac{d^n}{dx^n} \left[\displaystyle\int_x^\infty f(y)(y-x)^{n-\nu-1} dy\right] = \lim_{\delta \downarrow 0} W_\delta^{\nu,n}[f](x) \qquad \mathrm{Re}\,\nu < n$

where

(3.4) $\qquad W_\delta^{\nu,n}[f](x) := (-1)^n D_\delta^n[W^{\nu-n}[f]](x) = (-1)^n W^{\nu-n}[D_\delta^n[f]](x)$

We call $W_\delta^{\nu,n}[f]$ the *approximate fractional orthogonal derivative*. The last equality in (3.4) follows from the fact that the operators $W^{\nu-n}$ and $D_\delta^n$ are convolution operators (i.e. operators which commute with translations) and therefore the integrals can be interchanged. Then there follows by substitution of (3.1b)

(3.5) $\qquad W_\delta^{\nu,n}[f](x) = \dfrac{(-1)^n k_n \Gamma(n+1)}{h_n \Gamma(n-\nu)} \dfrac{1}{\delta^n} \displaystyle\int_{-\infty}^{\infty} \left[\int_{x+\delta u}^{\infty} f(y)(y - x - \delta u)^{n-\nu-1} dy\right] p_n(u) d\mu(u) =$

$\qquad\qquad = \dfrac{(-1)^n k_n \Gamma(n+1)}{h_n \Gamma(n-\nu)} \dfrac{1}{\delta^n} \displaystyle\int_{-\infty}^{\infty} \left[\int_x^{\infty} f(y + \delta u)(y - x)^{n-\nu-1} dy\right] p_n(u) d\mu(u)$

With $y = x + \delta(s - u)$ there follows

(3.6) $\qquad W_\delta^{\nu,n}[f](x) = \dfrac{(-1)^n k_n \Gamma(n+1)}{h_n \Gamma(n-\nu)} \dfrac{1}{\delta^\nu} \displaystyle\int_{-\infty}^{\infty} \left[\int_u^{\infty} f(x + \delta s)(s - u)^{n-\nu-1} ds\right] p_n(u) d\mu(u)$

Now we consider three typical cases for the orthogonality interval:
- Finite orthogonality interval $[-1, 1]$.
- Infinite orthogonality interval $[0, \infty)$.
- Infinite orthogonality interval $(-\infty, \infty)$.

- Finite orthogonality interval $[-1, 1]$.

The integral (3.6) can be split up in two integrals as

(3.7) $\qquad W_\delta^{\nu,n}[f](x) = \dfrac{(-1)^n k_n \Gamma(n+1)}{h_n \Gamma(n-\nu)} \dfrac{1}{\delta^\nu} \displaystyle\int_{-1}^{1} \left[\int_u^{1} f(x + \delta y)(y - u)^{n-\nu-1} dy\right] p_n(u) d\mu(u) +$

$\qquad\qquad + \dfrac{(-1)^n k_n \Gamma(n+1)}{h_n \Gamma(n-\nu)} \dfrac{1}{\delta^\nu} \displaystyle\int_{-1}^{1} \left[\int_1^{\infty} f(x + \delta y)(y - u)^{n-\nu-1} dy\right] p_n(u) d\mu(u)$

Interchanging the double integrals in both terms of (3.7) gives

$$(3.8) \quad W_\delta^{v,n}[f](x) = \frac{(-1)^n k_n \Gamma(n+1)}{h_n \Gamma(n-v)} \frac{1}{\delta^v} \int_{-1}^{1} f(x+\delta y)\left[\int_{-1}^{y} p_n(u)(y-u)^{n-v-1} d\mu(u)\right] dy +$$
$$+ \frac{(-1)^n k_n \Gamma(n+1)}{h_n \Gamma(n-v)} \frac{1}{\delta^v} \int_{1}^{\infty} f(x+\delta y)\left[\int_{-1}^{1} p_n(u)(y-u)^{n-v-1} d\mu(u)\right] dy$$

- Infinite orthogonality interval $[0,\infty)$.

The outer integral in (3.6) is now from $0$ to $\infty$. Interchanging the double integral gives

$$(3.9) \quad W_\delta^{v,n}[f](x) = \frac{(-1)^n k_n \Gamma(n+1)}{h_n \Gamma(n-v)} \frac{1}{\delta^v} \int_{0}^{\infty} f(x+\delta y)\left[\int_{0}^{y} (y-u)^{n-v-1} p_n(u) d\mu(u)\right] dy$$

- Infinite orthogonality interval $(-\infty,\infty)$.

The outer integral in (3.6) is now from $-\infty$ to $\infty$. Interchanging the double integral gives

$$(3.10) \quad W_\delta^{v,n}[f](x) = \frac{(-1)^n k_n \Gamma(n+1)}{h_n \Gamma(n-v)} \frac{1}{\delta^v} \int_{-\infty}^{\infty} f(x+\delta y)\left[\int_{-\infty}^{y} (y-u)^{n-v-1} p_n(u) d\mu(u)\right] dy$$

In formulas (3.8), (3.9) and (3.10) the integrals inside the square brackets can be expected to be computable analytically or numerically for concrete measures $\mu$ and orthogonality polynomials $p_n$.

Remark 3.2. We have derived these formulas for all values of $v$ with $\operatorname{Re} v < n$, but our main usage of them will be for $\operatorname{Re} v \geq 0$, where they are approximate fractional derivatives.

Remark 3.3. In this paper we assume for simplicity that the orthogonal derivative is equal to the ordinary derivative, i.e. that both sides of (3.1a) are well defined and equal, which is certainly true under the assumptions of Theorem 3.1. It should be noted that the definition of the orthogonal derivative is valid for a wider class of functions then the ordinary derivative. So if the limit of the right-hand side of (3.1a) exists and $f^n(x)$ does not exist, we still call the right-hand side of (3.1a) the orthogonal derivative. For example the ordinary derivative does not exist for the function $f(x) = |x|$ for $x = 0$, but the orthogonal derivative for this function does exist. For less trivial examples see Diekema and Koornwinder (2012), section 3.8.

Remark 3.4. Instead of the Weyl integral we could have worked with the Riemann-Liouville integral. For example the following formulas can be obtained

$$(3.11) \quad R_\delta^{v,n}[f](x) = R^{v-n}\left[D_\delta^n[f]\right](x) \qquad R^v[f](x) = \lim_{\delta \downarrow 0} R_\delta^{v,n}[f](x)$$

For the orthogonality interval $[-1,1]$ there follows

$$(3.12) \quad R_\delta^{v,n}[f](x) = \frac{k_n \Gamma(n+1)}{h_n \Gamma(n-v)} \frac{1}{\delta^v} \int_{-1}^{1} f(x+\delta y)\left[\int_{y}^{1} (u-y)^{n-v-1} p_n(u) d\mu(u)\right] dy +$$
$$+ \frac{k_n \Gamma(n+1)}{h_n \Gamma(n-v)} \frac{1}{\delta^v} \int_{-\infty}^{-1} f(x+\delta y)\left[\int_{-1}^{1} (u-y)^{n-v-1} p_n(u) d\mu(u)\right] dy$$

For the orthogonality interval $(-\infty,\infty)$ there follows

$$(3.13) \quad R_\delta^{v,n}[f](x) = \frac{k_n \Gamma(n+1)}{h_n \Gamma(n-v)} \frac{1}{\delta^v} \int_{-\infty}^{\infty} f(x+\delta y)\left[\int_{y}^{\infty} (u-y)^{n-v-1} p_n(u) d\mu(u)\right] dy$$

All results in this paper could also have been equivalently formulated in terms of $R_\delta^{v,n}$.

## 4 The fractional Weyl transform for the Jacobi derivative.

In this section we apply formula (3.8) to the Jacobi polynomials, where

(4.1a) $\quad p_n(x) = P_n^{(\alpha,\beta)}(x)$

(4.1b) $\quad d\mu(x) = w(x)dx$ with $w(x) = (1-x)^\alpha (1+x)^\beta \qquad \alpha > -1 \qquad \beta > -1$

(4.1c) $\quad \dfrac{h_n}{k_n} = 2^{n+\alpha+\beta+1} \dfrac{\Gamma(n+\alpha+1)\Gamma(n+\beta+1)}{\Gamma(2n+\alpha+\beta+2)}$

(4.1d) $\quad P_n^{(\alpha,\beta)}(-x) = (-1)^n P_n^{(\beta,\alpha)}(x)$

For the resulting transform, which we call the *approximate fractional Jacobi derivative*, we write $W_\delta^{v,n} = W_{\delta,\alpha,\beta}^{v,n}$. Before computing this transform we observe the following proposition

Proposition 4.1. For the approximate fractional Jacobi derivative the following formula is valid

(4.2) $\quad W_{\delta,\alpha,\beta}^{v,n}[f](x) = W_{\delta,\alpha+1,\beta+1}^{v,n-1}[f](x)$

These iteration goes on so long as $v < n - 1$.

Proof: From Diekema and Koornwinder (2012) formula (3.10) it follows that

(4.3) $\quad D_{\delta,\alpha,\beta}^n[f](x) = \dfrac{\int_{-1}^{1} f^{(n)}(x+\delta\xi)(1-\xi)^{n+\alpha}(1+\xi)^{n+\beta}d\xi}{\int_{-1}^{1}(1-\xi)^{n+\alpha}(1+\xi)^{n+\beta}d\xi} = D_{\delta,\alpha+1,\beta+1}^{n-1}[f'](x)$

From (2.8b) we obtain straight forward

(4.4) $\quad W_{\delta,\alpha,\beta}^{v,n}[f](x) = (-1)^n D_{\delta,\alpha,\beta}^n\left[W^{v-n}[f]\right](x) = (-1)^n D_{\delta,\alpha,\beta}^{n-1}\left(W^{v-n}[f]\right)'(x) =$

$\quad = (-1)^{n-1} D_{\delta,\alpha+1,\beta+1}^{n-1}\left[W^{v-n+1}[f]\right](x) = W_{\delta,\alpha+1,\beta+1}^{v,n-1}[f](x) \quad \square$

Remark 4.2. For the Laguerre and Hermite derivatives there are identical formulas. The formula for the Hermite derivative has no parameters.

For the computation of $W_{\delta,\alpha,\beta}^{v,n}$ we substitute (4.1) in (3.8).

(4.5) $\quad W_{\delta,\alpha,\beta}^{v,n}[f](x) = \dfrac{(-1)^n \Gamma(2n+\alpha+\beta+2)}{2^{n+\alpha+\beta+1}\Gamma(n+\alpha+1)\Gamma(n+\beta+1)} \dfrac{\Gamma(n+1)}{\Gamma(n-v)} \dfrac{1}{\delta^v}\left[I_1(x) + I_2(x)\right]$

where

(4.6a) $\quad I_1(x) = \int_{-1}^{1} f(x+\delta y) J_1(y) dy$

(4.6b) $\quad I_2(x) = \int_{1}^{\infty} f(x+\delta y) J_2(y) dy$

with

(4.7a) $\quad J_1(y) = \int_{-1}^{y} P_n^{(\alpha,\beta)}(u)(1-u)^\alpha(1+u)^\beta(y-u)^{n-v-1}du \quad -1 \leq y \leq 1$

(4.7b) $\quad J_2(y) = \int_{-1}^{1} P_n^{(\alpha,\beta)}(u)(1-u)^\alpha(1+u)^\beta(y-u)^{n-v-1}du \quad 1 \leq y \leq \infty$

For the computation of $J_1(y)$ the Rodrigues formula for Jacobi polynomials is used:

(4.8) $\quad P_n^{(\alpha,\beta)}(u) = \dfrac{1}{(-1)^n 2^n n!(1-u)^\alpha(1+u)^\beta} \dfrac{d^n}{du^n}\left[(1-u)^{n+\alpha}(1+u)^{n+\beta}\right]$

Substitution in (4.7a) yields:

(4.9) $\quad J_1(y) = \dfrac{(-1)^n}{2^n n!} \int_{-1}^{y} \dfrac{d^n}{du^n}\left[(1-u)^{n+\alpha}(1+u)^{n+\beta}\right](y-u)^{n-v-1}du$

Repeated integration by parts ($n$ times) gives:

(4.10) $\quad J_1(y) = \dfrac{(-1)^n}{2^n n!} \dfrac{\Gamma(n-v)}{\Gamma(-v)} \int_{-1}^{y}(1-u)^{n+\alpha}(1+u)^{n+\beta}(y-u)^{-v-1}du$

Substitution of the variable $u = (1+y)w - 1$, using Andrews et al. (1999) Theorem 2.2.1 and Olver et al. (2010) formula 15.8.1. in the third identity, gives:

(4.11) $\quad J_1(y) = (-1)^n 2^\alpha (1+y)^{n+\beta-v} \dfrac{\Gamma(n+\beta+1)\Gamma(n-v)}{n!\Gamma(n-v+\beta+1)} F\left(\begin{array}{c} n+\beta+1,-\alpha-n \\ n-v+\beta+1 \end{array}; \dfrac{1+y}{2}\right) =$

$$= (-1)^n \frac{\Gamma(n+\beta+1)\Gamma(n-\nu)}{2^{n-\nu}n!\Gamma(n-\nu+\beta+1)}(1-y)^{n+\alpha-\nu}(1+y)^{n+\beta-\nu}F\left(\begin{array}{c}-\nu, 2n-\nu+\alpha+\beta+1\\ n-\nu+\beta+1\end{array}; \frac{1+y}{2}\right)$$

where $F$ is the Gauss hypergeometric function.

Similarly the integral for $J_2(y)$ can be computed. The result is

$$(4.12) \quad J_2(y) = (-1)^n \frac{2^{n+\alpha+\beta+1}}{(y+1)^{\nu+1}} \frac{\Gamma(n-\nu)}{\Gamma(-\nu)n!} \frac{\Gamma(n+\alpha+1)\Gamma(n+\beta+1)}{\Gamma(2n+\alpha+\beta+2)} F\left(\begin{array}{c}\nu+1, n+\beta+1\\ 2n+\alpha+\beta+2\end{array}; \frac{2}{y+1}\right)$$

Substitution of $I_1(y)$ and $I_2(y)$ yields our result as the following theorem

**Theorem 4.3.** For the approximate fractional Jacobi derivative the following formula is valid:

$$(4.13) \quad W_{\delta,\alpha,\beta}^{\nu,n}[f](x) = \frac{1}{\Gamma(-\nu)} \frac{1}{\delta^\nu} \int_1^\infty f(x+\delta y) \frac{1}{(y+1)^{\nu+1}} F\left(\begin{array}{c}\nu+1, n+\beta+1\\ 2n+\alpha+\beta+2\end{array}; \frac{2}{y+1}\right) dy +$$

$$+ \frac{(-1)^n \Gamma(2n+\alpha+\beta+2)}{2^{2n-\nu+\alpha+\beta+1}\Gamma(n+\alpha+1)\Gamma(n-\nu+\beta+1)} \frac{1}{\delta^\nu}$$

$$\int_{-1}^1 f(x+\delta y)(1-y)^{n+\alpha-\nu}(1+y)^{n+\beta-\nu} F\left(\begin{array}{c}-\nu, 2n-\nu+\alpha+\beta+1\\ n-\nu+\beta+1\end{array}; \frac{1+y}{2}\right) dy$$

with $f(x) = O(\mathrm{Re}(x^{\nu-n-\varepsilon}))$ as $x \to \infty$ and $\mathrm{Re}\,\nu < n$. By continuity formula (4.13) is also valid as $\mathrm{Re}\,\nu = n$.

From this formula one can see once more the validity of Proposition 4.1. For $\nu = 0, 1, 2, \ldots, n$ the hypergeometric function in the first integral can be written as a Jacobi polynomial. The first term vanishes, because of the Gamma function $\Gamma(-\nu) = \Gamma(-n)$. In the second term the hypergeometric function can be written as $P_n^{(\beta,\alpha)}(-y)$. Using (4.1d) and taking the limit for $\delta \downarrow 0$ there remains the orthogonal derivative associated with the Jacobi polynomials

$$(4.14) \quad D^n[f](x) = \frac{\Gamma(2n+\alpha+\beta+2)\Gamma(n+1)}{2^{n+\alpha+\beta+1}\Gamma(n+\alpha+1)\Gamma(n+\beta+1)} \lim_{\delta \downarrow 0} \frac{1}{\delta^n} \int_{-1}^1 f(x+\delta y)(1-y)^\alpha(1+y)^\beta P_n^{(\alpha,\beta)}(y) dy$$

which is exactly (3.1b) for the case of the Jacobi polynomials.

For the special case that $\alpha \to \alpha - 1/2$ and $\beta \to \alpha - 1/2$ the Jacobi polynomials become the Gegenbauer polynomials. In that case the hypergeometric functions in (4.13) can be expressed in terms of associated Legendre functions. For the hypergeometric function in the first term of (4.13) Olver et al. (2010) formula 15.8.13 together with formula 14.3.7 gives:

$$(4.15) \quad \frac{1}{(y+1)^{\nu+1}} F\left(\begin{array}{c}\nu+1, n+\alpha+\frac{1}{2}\\ 2n+2\alpha+1\end{array}; \frac{2}{y+1}\right) = (-1)^n \frac{2^{n+\alpha+1/2}\Gamma(n+\alpha+1)}{e^{i\pi(1/2-\alpha+\nu)}\sqrt{\pi}\,\Gamma(\nu+1)} \frac{Q_{n+\alpha-1/2}^{1/2-\alpha-n+\nu}(y)}{\left(\sqrt{y^2-1}\right)^{1/2-\alpha-n+\nu}}$$

$Q_\nu^\mu(y)$ is the associated Legendre of the second kind. For the hypergeometric function in the second term of (4.13) Olver et al. (2010) formula 15.8.1 together with formula 14.3.1 gives:

$$(4.16) \quad F\left(\begin{array}{c}-\nu, 2n-\nu+2\alpha\\ n-\nu+\alpha+\frac{1}{2}\end{array}; \frac{1-y}{2}\right) = \Gamma\left(n-\nu+\alpha+\frac{1}{2}\right) 2^{n+\alpha-\nu-1/2} \frac{P_{n+\alpha-1/2}^{1/2-\alpha-n+\nu}(y)}{\left(\sqrt{1-y^2}\right)^{n+\alpha-\nu-1/2}}$$

$P_\nu^\mu(y)$ is the associated Legendre function of the first kind on the cut. Substitution gives for the *approximate fractional Gegenbauer* derivative:

$$(4.17) \quad W_{\delta,\alpha}^{\nu,n}[f](x)_G = (-1)^n \frac{\Gamma(n+\alpha+1)2^{n+\alpha-1/2}}{\sqrt{\pi}} \frac{1}{\delta^\nu} \int_{-1}^1 f(x-\delta y) \frac{P_{n+\alpha-1/2}^{1/2-\alpha-n+\nu}(y)}{\left(\sqrt{1-y^2}\right)^{1/2-\alpha-n+\nu}} dy -$$

$$- (-1)^n \frac{e^{i\pi(\alpha-1/2)}2^{n+\alpha+1/2}\Gamma(n+\alpha+1)\sin(\nu\pi)}{\pi\sqrt{\pi}} \frac{1}{\delta^\nu} \int_1^\infty f(x-\delta y) \frac{Q_{n+\alpha-1/2}^{1/2-\alpha-n+\nu}(y)}{\left(\sqrt{y^2-1}\right)^{1/2-\alpha-n+\nu}} dy \qquad \mathrm{Re}\,\nu \le n$$

Further specialization to $\alpha = 1/2$ gives the *approximate fractional Legendre* derivative:

$$\text{(4.18)} \quad W_{\delta,1/2}^{v,n}[f](x)_P = (-1)^n \frac{2^n}{\sqrt{\pi}} \Gamma\left(n + \frac{3}{2}\right) \frac{1}{\delta^v} \int_{-1}^{1} f(x - \delta y) \frac{P_n^{v-n}(y)}{\left(\sqrt{1-y^2}\right)^{v-n}} dy -$$

$$- (-1)^n \frac{2^{n+1} \sin(v\pi)}{\pi\sqrt{\pi}} \Gamma\left(n + \frac{3}{2}\right) \frac{1}{\delta^v} \int_{1}^{\infty} f(x - \delta y) \frac{Q_n^{v-n}(y)}{\left(\sqrt{y^2-1}\right)^{v-n}} dy \qquad \text{Re}\, v \le n$$

Remark 4.4. For the approximate fractional Gegenbauer derivative following the Riemann-Liouville definition it can be shown that with $\text{Re}\, v \le n$

$$\text{(4.19)} \quad R_{\delta,\alpha}^{v,n}[f](x)_G = \frac{2^{n+\alpha-1/2}\Gamma(n+\alpha+1)}{\sqrt{\pi}} \frac{1}{\delta^v} \int_{-1}^{1} f(x+\delta y) \frac{P_{n+\alpha-1/2}^{1/2-\alpha-n+v}(y)}{\left(\sqrt{1-y^2}\right)^{1/2-\alpha-n+v}} dy -$$

$$- (-1)^n \frac{e^{i\pi(\alpha-1/2)} 2^{n+\alpha+1/2}\Gamma(n+\alpha+1)\sin(v\pi)}{\pi\sqrt{\pi}} \frac{1}{\delta^v} \int_{1}^{\infty} f(x+\delta y) \frac{Q_{n+\alpha-1/2}^{1/2-\alpha-n+v}(y)}{\left(\sqrt{y^2-1}\right)^{1/2-\alpha-n+v}} dy$$

Remark 4.5. In the same way as we did in the Jacobi case one can calculate the *approximate fractional Laguerre derivative*. Then we start with formula (3.9). There follows

$$\text{(4.20)} \quad W_{\delta,\alpha}^{v,n}[f](x)_L = \frac{1}{\Gamma(n-v+\alpha+1)} \frac{1}{\delta^v} \int_{0}^{\infty} f(x+\delta y) y^{n-v+\alpha} e^{-y} M(-v; n-v+\alpha+1; y) dy$$

where $M$ is the confluent hypergeometric function and $f(x) = O(\text{Re}(x^{v-n-\varepsilon}))$ as $x \to \infty$ with $\text{Re}\, v \le n$.

## 5  The transfer function for the approximate fractional orthogonal derivative.

In section 2 Theorem 2.1. we mentioned the Fourier transform of the fractional derivative. It is possible to extend that theorem to the approximate fractional orthogonal derivative.

As in Theorem 2.1 let $F(\omega)$ and $G(\omega)$ be the Fourier transforms of $f$ and $g$ respectively. From (3.1b), with $f$ absolutely integrable on $\mathbb{R}$, and $g = D_\delta^n[f]$, we immediately derive for the transfer function $H_\delta^n(\omega)$

$$\text{(5.1)} \quad H_\delta^n(\omega) = \frac{G(\omega)}{F(\omega)} = \frac{k_n n!}{h_n \delta^n} \int_{\mathbb{R}} p_n(\xi) e^{-i\omega\delta\xi} d\mu(\xi)$$

The exponential can be written as a sum. Because of the orthogonality property there follows

$$\text{(5.2)} \quad H_\delta^n(\omega) = \frac{k_n n!}{h_n \delta^n} \int_{\mathbb{R}} p_n(\xi) \sum_{k=n}^{\infty} \frac{(-i\omega\delta\xi)^k}{k!} d\mu(\xi) = \frac{k_n n!}{h_n} (-i\omega)^n \int_{\mathbb{R}} p_n(\xi) \xi^n \sum_{k=0}^{\infty} \frac{(-i\omega\delta\xi)^k}{(k+n)!} d\mu(\xi)$$

Now there follows the next theorem

Theorem 5.1. Let $F(\omega)$ and $G(\omega)$ be the Fourier transforms of the function $f$ and the function $g = W_\delta^{v,n}[f]$, respectively. Then (recall the convention after (2.11))

$$\text{(5.3)} \quad H_\delta^v(\omega) = \frac{k_n n!}{h_n \delta^n} (-1)^n (i\omega)^{v-n} \int_{\mathbb{R}} p_n(\xi) e^{-i\omega\delta\xi} d\mu(\xi) = \frac{k_n n!}{h_n} (i\omega)^v \int_{\mathbb{R}} p_n(\xi) \xi^n \sum_{k=0}^{\infty} \frac{(-i\omega\delta\xi)^k}{(k+n)!} d\mu(\xi)$$

Proof: Combine (5.1) and (5.2) with (3.4) and Theorem 2.1. and the result follows immediately.□

Since $\int_{\mathbb{R}} p_n(\xi) e^{-i\omega\delta\xi} d\mu(\xi)$ is explicitly known for the classical orthogonal polynomials, we get explicit formulas for the transfer function in these cases. Also for the case of discrete Hahn polynomials it is possible to compute the transfer function. We now treat two examples. First we compute the transfer function in the case of the Jacobi polynomials. Then we compute the transfer function in the case of the discrete Hahn polynomials.

1. The transfer function for the Jacobi polynomials.

For the Jacobi polynomials the following integral is known (Olver et al. (2010), 18.17.16)

$$\text{(5.4)} \quad \int_{-1}^{1} (1-x)^\alpha (1+x)^\beta P_n^{(\alpha,\beta)}(x) e^{ixy} dx =$$

$$= \frac{(iy)^n e^{iy}}{n!} 2^{n+\alpha+\beta+1} B(n+\alpha+1, n+\beta+1) M(n+\alpha+1, 2n+\alpha+\beta+2; -2iy)$$

Using (4.1b) and (4.1c) and substitution of $y = -\omega\delta$ yields for (5.3) after some calculation

(5.5) $\quad H_\delta^\nu(\omega) = (i\omega)^\nu e^{-i\omega\delta} M(n + \alpha + 1, 2n + \alpha + \beta + 2; 2i\omega\delta)$

Using (Olver et al. (2010), 13.7.2) for $\omega \to \infty$ there remains

(5.6) $\quad \left|H_\delta^\nu(\omega)\right| \sim const \cdot \omega^{\nu-n-\min(\alpha,\beta)}$

Hence. since $\nu < n$

(5.7) $\quad \lim_{\omega \to \infty} \left|H_\delta^\nu(\omega)\right| = 0$

For the Gegenbauer polynomials let $\alpha \to \alpha - \frac{1}{2}$ and $\beta \to \alpha - \frac{1}{2}$. Then

(5.8) $\quad H_\delta^\nu(\omega) = (i\omega)^\nu e^{-i\omega\delta} M\left(n + \alpha + \frac{1}{2}; 2n + 2\alpha + 1; 2i\omega\delta\right)$

The confluent hypergeometric function can be transformed into a Bessel function. Olver et al. (2010) formula 13.6.9 gives:

(5.9) $\quad M\left(\nu + \frac{1}{2}, 2\nu + 1; 2iz\right) = \Gamma(1 + \nu) e^{iz} \left(\frac{2}{z}\right)^\nu J_\nu(z)$

Substitution in (5.5) gives

(5.10) $\quad H_\delta^\nu(\omega) = (i\omega)^\nu \Gamma(n + \alpha + 1) \left(\frac{2}{\omega\delta}\right)^{n+\alpha} J_{n+\alpha}(\omega\delta)$

The Bessel function can be written as a summation (Olver et al. (2010), 10.2.2). After substitution there remains

(5.11) $\quad H_\delta^\nu(\omega) = (i\omega)^\nu \Gamma(n + \alpha + 1) \sum_{k=0}^{\infty} \frac{(-1)^k}{\Gamma(k + n + \alpha + 1) k!} \left(\frac{\omega\delta}{2}\right)^{2k}$

The right-hand side is real except for the factor $(i\omega)^\nu$. The variable $n$ is equal to the next integer greater then $\nu$. For a certain value of $\nu$ with $n - 1 \leq \nu \leq n$ the summation is independent of $\nu$. For the modulus of the transfer function we obtain

(5.12) $\quad \left|H_\delta^\nu(\omega)\right| = \omega^\nu \Gamma(n + \alpha + 1) \left|\sum_{k=0}^{\infty} \frac{(-1)^k}{\Gamma(k + n + \alpha + 1) k!} \left(\frac{\omega\delta}{2}\right)^{2k}\right|$

$\alpha$ is a free parameter. For $\alpha = \frac{1}{2}$ we get the transfer function of the approximate fractional Legendre derivative. Then (5.13) becomes

(5.13) $\quad \left|H_\delta^\nu(\omega)\right| = \omega^\nu \frac{\Gamma(2n + 2)}{2^n \Gamma(n + 1)} \frac{1}{(\omega\delta)^n} \left|j_n(\omega\delta)\right|$

where $j_n$ is a spherical Bessel function (Olver et al. (2010), 10.49.2). Special cases of $j_n$ are

(5.14a) $\quad j_0(z) = \frac{1}{z} \sin z$

(5.14b) $\quad j_1(z) = \frac{1}{z^2} (\sin z - z \cos z)$

(5.14c) $\quad j_2(z) = \frac{1}{z^3} \left((3 - z^2) \sin z - 3z \cos z\right)$

(5.14d) $\quad j_3(z) = \frac{1}{z^4} \left((15 - 6z^2) \sin z - (15 - z^2) z \cos z\right)$

The above formulas remains valid for $\nu = n$, where we get the transfer functions as in Diekema et al (2012), section 5.

Remark 5.2. From Theorem 5.1. we get

$$\text{(5.15)} \quad \frac{1}{\sqrt{2\pi}} \int_{\mathbb{R}} W^{v,n}_\delta[f](x)e^{i\omega x}dx = H^v_\delta(\omega)F(\omega)$$

with $H^v_\delta(\omega)$ explicitly and quickly computable in special cases, for instance the Jacobi case. Thus if

$$\text{(5.16)} \quad H^v_\delta(\omega) = \frac{1}{\sqrt{2\pi}} \frac{1}{\delta^v} \int_{\mathbb{R}} h^v_\delta(x)e^{i\sigma\omega x}dx$$

for some function $h$ then after the inverse Fourier transform there remains the convolution integral

$$\text{(5.17)} \quad W^{v,n}_\delta[f](x) = \frac{1}{\sqrt{2\pi}} \frac{1}{\delta^v} \int_{\mathbb{R}} f(x - \delta y)h^v_\delta(y)dy$$

If we compare, in the case of a finite orthogonality, this last formula with (3.8) then we can read off $h$ from (3.8). In particular, in the Jacobi case we can read off $h$ from (4.13). This gives

$$\text{(5.18)} \quad h^v_\delta(y) = \begin{cases} \dfrac{\sqrt{2\pi}}{\Gamma(-v)\delta^{v+1}} \dfrac{1}{(1+y)^{v+1}} F\!\left(\begin{matrix} v+1, n+\beta+1 \\ 2n+\alpha+\beta+2 \end{matrix}; \dfrac{2}{1+y}\right) & y > 1 \\[1em] \dfrac{\sqrt{2\pi}\,\Gamma(2n+\alpha+\beta+2)}{2^{2n-v+\alpha+\beta+1}\Gamma(n+\alpha+1)\Gamma(n-v+\beta+1)\delta^{v+1}} \\[0.5em] \quad (1+y)^{n+\beta-v}(1-y)^{n+\alpha-v} F\!\left(\begin{matrix} -v, 2n-v+\alpha+\beta+1 \\ n-v+\beta+1 \end{matrix}; \dfrac{1+y}{2}\right) & -1 < y < 1 \\[1em] 0 & y < -1 \end{cases}$$

Its Fourier transform $H^v_\delta(\omega)$ is given by (5.5).

2. The transfer function for the Hahn polynomials.

We start with formula (5.3) for the approximate fractional orthogonal derivative and apply this formula for the discrete case with weights $w(x)$ on points $x \in \{0, 1, \ldots, N\}$. This gives

$$\text{(5.19)} \quad H^v_\delta(\omega) = \frac{k_n n!}{h_n \delta^n}(-1)^n(i\omega)^{v-n} \sum_{x=0}^{N} p_n(x)w(x)e^{-ix\delta\omega}$$

When applying the summation to the discrete Hahn polynomials we get

$$\text{(5.20)} \quad S = \sum_{x=0}^{N} p_n(x)w(x)e^{-ix\delta\omega} = \sum_{x=0}^{N} Q_n(x;\alpha,\beta,N)w(x)e^{-i\delta\omega x}$$

For the Hahn polynomials we use formula (9.5.10) of Koekoek et al. (2010)

$$\text{(5.21)} \quad Q_n(x;\alpha,\beta,N)w(x) = \frac{(-1)^n(\beta+1)_n}{(-N)_n} \nabla^n_x\!\left(\frac{(\alpha+n+1)_x}{x!} \frac{(\beta+n+1)_{N-n-x}}{(N-n-x)!}\right)$$

where

$$\text{(5.22)} \quad \nabla_x(f(x)) := f(x) - f(x-1)$$

and get

$$\text{(5.23)} \quad S = \frac{(-1)^n(\beta+1)_n}{(-N)_n} \sum_{x=0}^{N} e^{-i\delta\omega x}\nabla^n\!\left(\frac{(\alpha+n+1)_x}{x!} \frac{(\beta+n+1)_{N-n-x}}{(N-n-x)!}\right)$$

To compute this sum we use partial summation. The general formula for partial summation is

$$\text{(5.24)} \quad \sum_{j=0}^{N}(\nabla f)(j)g(j) = f(N)g(N) - f(-1)g(0) - \sum_{j=0}^{N-1} f(j)(\nabla g)(j+1)$$

In particular, if $f(N) = f(-1) = 0$ then

(5.25) $$\sum_{j=0}^{N}(\nabla f)(j)g(j) = -\sum_{j=0}^{N-1}f(j)(\nabla g)(j+1)$$

Applying this formula $n$ times and remarking that $f(N-n) = f(-1) = 0$ gives

(5.26) $$\sum_{j=0}^{N}(\nabla f)^n(j)g(j) = (-1)^n\sum_{j=0}^{N-n}f(j)(\nabla g)^n(j+n)$$

Applying this formula to (5.23) gives

(5.27) $$S = \frac{(\beta+1)_n}{(-N)_n}\sum_{x=0}^{N-n}\nabla^n(e^{-i\delta\omega(x+n)})\frac{(\alpha+n+1)_x}{x!}\frac{(\beta+n+1)_{N-n-x}}{(N-n-x)!}$$

The iterated differences of the exponential are known. This gives

(5.28) $$S = \frac{(\beta+1)_n}{(-N)_n}(e^{-i\delta\omega}-1)^n\sum_{x=0}^{N-n}e^{-i\delta\omega x}\frac{(\alpha+n+1)_x}{x!}\frac{(\beta+n+1)_{N-n-x}}{(N-n-x)!} =$$
$$= \frac{(\beta+1)_n}{(-N)_n}(e^{-i\delta\omega}-1)^n\sum_{x=0}^{N-n}\frac{(\alpha+n+1)_x(\beta+n+1)_{N-n}}{(-\beta-N)_x}\frac{(-N+n)_x}{(N-n)!}\frac{(e^{-i\delta\omega})^x}{x!} =$$
$$= \frac{(\beta+1)_n(\beta+n+1)_{N-n}}{(-N)_n(N-n)!}(e^{-i\delta\omega}-1)^n\sum_{x=0}^{N-n}\frac{(\alpha+n+1)_x(-N+n)_x}{(-\beta-N)_x}\frac{(e^{-i\delta\omega})^x}{x!} =$$
$$= \frac{(\beta+1)_N}{N!}(1-e^{-i\delta\omega})^n F\begin{pmatrix}-N+n,\alpha+n+1\\-\beta-N\end{pmatrix};e^{-i\delta\omega}$$

Then for the transfer function there remains

(5.29) $$H^\gamma_\delta(\omega) = \frac{k_n n!}{h_n \delta^n}(-1)^n(i\omega)^{\nu-n}\frac{(\beta+1)_N}{N!}(1-e^{-i\delta\omega})^n F\begin{pmatrix}-N+n,\alpha+n+1\\-\beta-N\end{pmatrix};e^{-i\delta\omega}$$

For the Hahn polynomials tables 18.19.1 and 18.19.2 in Olver et al (2010) give

(5.30) $$\frac{k_n}{h_n} = (-1)^n\frac{(2n+\alpha+\beta+1)}{(\beta+1)_n}\frac{(n+\alpha+\beta+1)_n}{(n+\alpha+\beta+1)_{N+1}}\frac{N!}{n!}$$

Substitution in (5.25) gives

(5.31) $$H^\gamma_\delta(\omega) = (i\omega)^{\nu-n}\frac{(1-e^{-i\delta\omega})^n}{\delta^n}\frac{\Gamma(N+\beta+1)\Gamma(2n+\alpha+\beta+2)}{\Gamma(n+\beta+1)\Gamma(N+n+\alpha+\beta+2)}F\begin{pmatrix}-N+n,\alpha+n+1\\-\beta-N\end{pmatrix};e^{-i\delta\omega}$$

Remark 5.3. From the summation formula (5.20) together with (5.28) one can derive the integral formula when letting $N \to \infty$. For this purpose we apply table 18.19 and the formulas (15.8.7) and (18.20.5) of Olver et al (2010) to (5.28). With $\delta = -1$ this gives

(5.32) $$\sum_{x=0}^{N}{}_3F_2\begin{pmatrix}-n,n+\alpha+\beta+1,-x\\\alpha+1,-N\end{pmatrix};1\frac{(\alpha+1)_x(\beta+1)_{N-x}}{x!(N-x)!}e^{i\omega x} =$$
$$= \frac{(\beta+1)_n(2n+\alpha+\beta+2)_{N-n}}{N!}(1-e^{i\omega})^n F\begin{pmatrix}-N+n,\alpha+n+1\\2n+\alpha+\beta+2\end{pmatrix};1-e^{i\omega}$$

The summation can also be done over the set $x \in \frac{1}{N}\{0, 1, \ldots, N\}$. The $N+1$ points then lie in the interval $[0, 1]$. This gives

(5.33) $$\sum_{x=0,1/N,..1}{}_3F_2\begin{pmatrix}-n,n+\alpha+\beta+1,-Nx\\;\alpha+1,-N\end{pmatrix};1\frac{(\alpha+1)_{Nx}(\beta+1)_{N(1-x)}}{(Nx)!(N-Nx)!}e^{i\omega x} =$$

$$= \frac{(\beta+1)_n (2n+\alpha+\beta+2)_{N-n}}{N!} (1-e^{i\omega/N})^n F\left(\begin{array}{c}-N+n, \alpha+n+1 \\ 2n+\alpha+\beta+2\end{array}; 1-e^{i\omega/N}\right)$$

This formula can be written as

(5.34) 
$$\frac{\Gamma(\alpha+1)\Gamma(\beta+1)}{N^{\alpha+\beta+1}} \sum_{x=0,1/N,\ldots 1} {}_3F_2\left(\begin{array}{c}-n, n+\alpha+\beta+1, -Nx \\ \alpha+1, -N\end{array}; 1\right) \frac{(\alpha+1)_{Nx}(\beta+1)_{N(1-x)}}{(Nx)!(N-Nx)!} e^{i\omega x} =$$
$$= \frac{\Gamma(\alpha+1)\Gamma(\beta+1)}{N^{n+\alpha+\beta+1}} \frac{(\beta+1)_n (2n+\alpha+\beta+2)_{N-n}}{N!} \left(N(1-e^{i\omega/N})\right)^n F\left(\begin{array}{c}-N+n, \alpha+n+1 \\ 2n+\alpha+\beta+2\end{array}; \frac{N(e^{i\omega/N}-1)}{-N}\right)$$

Taking the limit for $N \to \infty$ and using (Olver et al. (2010), 5.11.12) gives formally

(5.35)
$$\int_0^1 F\left(\begin{array}{c}-n, n+\alpha+\beta+1 \\ \alpha+1\end{array}; x\right) x^\alpha (1-x)^\beta e^{i\omega x} dx =$$
$$= \frac{\Gamma(\alpha+1)\Gamma(n+\beta+1)}{\Gamma(2n+\alpha+\beta+2)} (-i\omega)^n M(n+\alpha+1; 2n+\alpha+\beta+2; i\omega)$$

This is the same formula as (5.4) after some transformation.

Remark 5.4. Tseng et al. (2000) compute the fractional derivative using the Fourier transform. They also design a fractional digital filter. They demonstrate the filter performance by using a block signal as input. In our opinion this demonstration is less convincing because by absence of an integrating factor the high frequencies are not suppressed, and because the plots are not given as log-log plots. One should compare with our approach in Section 11 where we give a log-log plot of the modulus of the transfer function.

## 6   Deriving a formula for a fractional derivative starting with a suitable transfer function.

For the second expression in (5.3) of the transfer function of the approximate fractional orthogonal derivative we will see that if $\delta\omega \to 0$, then

(6.1)   $$H^\nu_\delta(\omega) = (i\omega + 0)^\nu \left(1 + O(\delta\omega)\right) = \frac{1}{\delta^\nu} (i\delta\omega + 0)^\nu \left(1 + O(\delta\omega)\right)$$

Conversely, if we have some explicit function $H^\nu_\delta(\omega) = (i\omega)^\nu \left(1 + O(\delta\omega)\right)$ as $\delta\omega \downarrow 0$ of which the inverse Fourier transform $h(x)$ exists and is explicitly known, then the convolution product of $f$ with $h$ gives a formula for an approximate fractional derivative. The next derivation demonstrates this method.
Set

(6.2)   $$F(\omega) = \frac{1}{\sqrt{2\pi}} \int_{-\infty}^{\infty} e^{i\omega x} f(x) dx$$

When taking the inverse Fourier transform of (2.11) there follows for the Weyl derivative

(6.3)   $$W^\nu[f](x) = \frac{1}{\sqrt{2\pi}} \int_{-\infty}^{\infty} e^{-i\omega x} (i\omega + 0)^\nu F(\omega) d\omega$$

Now we look for a function with the following property

(6.4)   $$\lim_{\delta \downarrow 0} H^\nu_\delta(i\omega + 0) = (i\omega + 0)^\nu$$

Substitution in (6.3) with $W^\nu[f](x) = \lim_{\delta \downarrow 0} W^\nu_\delta[f](x)$ yields (care should be taken with the interchanging of the limit and the integral)

(6.5)   $$W^\nu_\delta[f](x) = \frac{1}{\sqrt{2\pi}} \int_{-\infty}^{\infty} e^{-i\omega x} H^\nu_\delta(i\omega + 0) F(\omega) d\omega$$

Suppose $H^\nu_\delta$ can be written as a Fourier transform

(6.6)   $$H^\nu_\delta(i\omega + 0) = \frac{1}{\delta^\nu} \frac{1}{\sqrt{2\pi}} \int_{-\infty}^{\infty} h^\nu_\delta(y) e^{i\delta\omega y} dy$$

Then for the approximate Weyl derivative there follows the convolution integral

$$(6.7) \quad W^{v}_{\delta}[f](x) = \frac{1}{\sqrt{2\pi}} \frac{1}{\delta^v} \int_{-\infty}^{\infty} h^{v}_{\delta}(y) f(x - \delta y) dy$$

Taking the limit for $\delta \downarrow 0$ yields

$$(6.8) \quad \frac{d^v}{dx^v}[f(x)] = \lim_{\delta \to 0} W^{v}_{\delta}[f](x) = \frac{1}{\sqrt{2\pi}} \lim_{\delta \to 0} \frac{1}{\delta^v} \int_{-\infty}^{\infty} h^{v}_{\delta}(y) f(x - \delta y) dy$$

The author started his original work on the fractional derivative with the special case $H(\omega) = \omega^{-\beta} J_a(\omega)$. Taking the inverse Fourier transform he derived the formula for the fractional Gegenbauer derivative (4.17) with two free parameters. No use was made of the orthogonality property of the Gegenbauer polynomials.

Because the Bessel function is a special case of the confluent hypergeometric function a next extension arises when using the inverse Fourier transform of the function $H(\omega) = \omega^{b-1} M(a, c; ik\omega)$. Then there are three free parameters. By (5.5) the approximate fractional Jacobi derivative will arise. A possible next extension is the use of the generalized confluent hypergeometric function. The transfer function can be written as $H(\omega) = \omega^{\gamma-1} {}_A F_A[(a), (b); ik\omega]$. There are $2A + 1$ free parameters. Possibly extensions could involve the Meijer $G$-function and the Fox $H$-function. See Kiryakova (1994).

The main properties of the function $H(\omega)$ are that its inverse Fourier transform exists and $H^{v}_{\delta}(\omega) = (i\omega + 0)^v \big(1 + O(\delta\omega)\big)$ as $\delta\omega \downarrow 0$.

As an example of the method we derive the formula for the fractional derivative using the inverse Fourier transform of the confluent hypergeometric function $H(\omega) = (i\omega + 0)^{b-1} M(a; c; i\omega)$.

For this transform there can be derived (see Appendix A)

$$(6.9) \quad \mathcal{F}^{-1}\left[(i\omega + 0)^{b-1} M(a, c; i\omega)\right](y) =$$

$$= \begin{cases} 0 & 1 < y \\ \sqrt{2\pi} \dfrac{\Gamma(c)}{\Gamma(a)\Gamma(1-a-b+c)} y^{a-b}(1-y)^{c-a-b} F\begin{pmatrix} 1-b, c-b \\ 1-a-b+c \end{pmatrix}; 1-y \end{pmatrix} & 0 < y < 1 \\ \dfrac{\sqrt{2\pi}}{\Gamma(1-b)} \dfrac{1}{(1-y)^b} F\begin{pmatrix} b, c-a \\ c \end{pmatrix}; \dfrac{1}{1-y} \end{pmatrix} & y < 0 \end{cases}$$

with conditions for the parameters: $0 < b < \min(a, c - a)$. Application of (6.7) with $v = b - 1$ gives

$$(6.10) \quad \frac{d^v}{dx^v}[f(x)] = \frac{1}{\Gamma(-v)} \lim_{\delta \to 0} \frac{1}{\delta^v} \int_{-\infty}^{0} \frac{1}{(1-y)^{v+1}} F\begin{pmatrix} v+1, c-a \\ c \end{pmatrix}; \frac{1}{1-y} \end{pmatrix} f(x - \delta y) dy +$$
$$+ \frac{\Gamma(c)}{\Gamma(a)\Gamma(c-a-v)} \lim_{\delta \to 0} \frac{1}{\delta^v} \int_{0}^{1} y^{a-v-1}(1-y)^{c-a-v-1} F\begin{pmatrix} -v, c-v-1 \\ c-a-v \end{pmatrix}; 1-y \end{pmatrix} f(x - \delta y) dy$$

Because the fractional derivative should be calculated, we set $v > 0$. This gives $b > 1$. In formula (6.10) there are two free parameters, namely $a$ and $c$. For the conditions we find: $0 < v < \min(a, c - a) - 1$.

Replacing in both integrals $y$ by $\dfrac{1-y}{2}$ and $\delta$ by $2\delta$ gives

$$(6.11) \quad \frac{d^v}{dx^v}[f(x)] = \frac{1}{\Gamma(-v)} \lim_{\delta \to 0} \frac{1}{\delta^v} \int_{1}^{\infty} \frac{1}{(1+y)^{v+1}} F\begin{pmatrix} v+1, c-a \\ c \end{pmatrix}; \frac{2}{1+y} \end{pmatrix} f(x + \delta y) dy +$$
$$+ \frac{\Gamma(c)}{2^{c-v-1}\Gamma(a)\Gamma(c-a-v)} \lim_{\delta \to 0} \frac{1}{\delta^v} \int_{-1}^{1} (1-y)^{a-v-1}(1+y)^{c-a-v-1} F\begin{pmatrix} -v, c-v-1 \\ c-a-v \end{pmatrix}; \frac{1+y}{2} \end{pmatrix} f(x + \delta y) dy$$

For $a = n + \alpha + 1$ and $c = 2n + \alpha + \beta + 2$ this formula is exactly the formula of the fractional Jacobi derivative (4.13) by taking the limit for $\delta \to 0$.

Remark 6.1. For the first method we use consecutively formulas (3.4), (3.8) which leads to (4.13) and to (5.5). As a side result of this first method we get (5.5) as the Fourier tramsform of (5.18), see Remark 5.2. The resulting inverde Fourier transform coincides after appropiate substitutions with (6.9).

Remark 6.2. With this method no use is made of the orthogonality property. So for extending this method with for example the generalized confluent hypergeometric function one cannot see directly what functions should be used with the first method to obtain the same result. Then there remains the question if in that case there are some orthogonal polynomials (these maybe yet unknown) that fulfill the first method.

Remark 6.3. One of the properties of the orthogonal derivative is that for high frequencies the modulus of the transfer function goes to zero. Our example satisfy this condition.

## 7 The fractional derivative arising from some particular functions.

In this section a demonstration is given of the second method for derivation of a formula of the fractional derivative based on some special functions. All the examples use an inverse Fourier transform in order to calculate a formula for the fractional derivative.

**Example 1**.
Consider the function

(7.1) $\quad H_a^v(\omega) = (i\omega + 0)^v e^{-a|\omega|} \qquad v > -1 \qquad a > 0$

Clearly if $a \downarrow 0$ then $H(\omega) \to (i\omega + 0)^v$. We calculate the inverse Fourier transform

(7.2) $\quad h_a^v(y) = \frac{1}{\sqrt{2\pi}} \int_{-\infty}^{\infty} (i\omega + 0)^v e^{-a|\omega|} e^{-i\omega y} d\omega = \frac{e^{i\pi v/2}}{\sqrt{2\pi}} \int_0^{\infty} \omega^v e^{-a\omega} e^{-i\omega y} d\omega + \frac{e^{-i\pi v/2}}{\sqrt{2\pi}} \int_0^{\infty} \omega^v e^{-a\omega} e^{i\omega y} d\omega =$

$\quad = \frac{\Gamma(v+1)e^{i\pi/2}}{\sqrt{2\pi}} \left[ \frac{1}{(y+ia)^{v+1}} - \frac{1}{(y-ia)^{v+1}} \right]$

Hence, with

(7.3) $\quad F(x) = \frac{1}{\sqrt{2\pi}} \int_{-\infty}^{\infty} f(\omega) e^{i\omega x} d\omega$

we obtain

(7.4) $\quad \frac{d^v}{dx^v} f(x) = \frac{1}{\sqrt{2\pi}} \lim_{a \to 0} \int_{-\infty}^{\infty} F(\omega)(i\omega + 0)^v e^{-a|\omega|} e^{-i\omega y} d\omega =$

$\quad = \frac{\Gamma(v+1)}{2\pi i} \lim_{a \to 0} \int_{-\infty}^{\infty} \left[ \frac{1}{(y-ia)^{v+1}} - \frac{1}{(y+ia)^{v+1}} \right] f(x+y) dy$

Care should be taken when $a \downarrow 0$. For $-\infty < y \le 0$ there is a cut, so we has to use a well-chosen contour to calculate the integral. In that case there remains

(7.5) $\quad R^v[f](x) = \frac{1}{\sqrt{2\pi}} \int_{-\infty}^{\infty} F(\omega)(i\omega + 0)^v e^{-i\omega y} d\omega = \frac{\Gamma(v+1)}{2\pi i} \int_C f(x+y) \frac{1}{y^{v+1}} dy$

Here $C$ is the contour with

(7.6) $\quad C = \{y - i0 \mid -\infty < y < -r\} \cup \{z \in \mathbb{C} \mid |z| = r\} \cup \{y + i0 \mid -r > y > -\infty\}$

This formula is well-known in the literature. See for example Lavoie et al. (1976) and Zavada (1998). If $v < 0$ then there are no convergence problems.
After application of formula (Olver et al. (2010), 15.4.9) there remains

(7.7) $\quad \frac{d^v}{dx^v} f(x) = \frac{\Gamma(v+2)}{\pi} \lim_{a \to 0} \frac{1}{a^{v-1}} \int_{-\infty}^{\infty} \frac{1}{y^{v+2}} F\left( \frac{v+2}{2}, \frac{v+2}{2} + \frac{1}{2}; \frac{3}{2}; -\frac{a^2}{y^2} \right) f(x+y) dy$

The quadratic argument in the hypergeometric function suggests a connection with the pseudo Jacobi polynomials. They are orthogonal with the orthogonality property (formula (9.9.2) of Koekoek et al. (2010))

(7.8) $\quad \int_{-\infty}^{\infty} P_m(x;0,N) P_n(x;0,N)(1+x^2)^{-N-1} dx = 0 \qquad n \ne m, \ n,m = 0,1,\ldots,N$

Formula (9.9.2) of Koekoek et al. (2010) gives

(7.9) $\quad P_{2n}(x;0,N) = (-2i)^{2n}(-N)_{2n}(2n-2N-1)_{2n}F\left(\begin{array}{c}-2n, 2n-2N-1\\-N\end{array};\frac{1-ix}{2}\right)$

With Olver et al. (2010) formulas 15.8.18 and 15.8.7 there remains

(7.10) $\quad P_{2n}(x;0,2n) = \frac{(-1)^n}{2n+1}F\left(\begin{array}{c}-n,-n-\frac{1}{2}\\\frac{1}{2}\end{array};-x^2\right)$

and with Olver et al. (2010) formula 15.4.7 there remains with a scaling factor $a$

(7.11) $\quad P_{2n}\left(\frac{x}{a};0,2n\right) = \frac{(-1)^n}{2(2n+1)}((1+iax)^{2n+1}-(1-iax)^{2n+1}) =$

$$= \frac{1}{2(2n+1)ia}(a^2+x^2)^{2n+1}\left(\frac{1}{(x-ia)^{2n+1}}-\frac{1}{(x+ia)^{2n+1}}\right)$$

Substituting $v = 2n$ in (7.4) gives

(7.12) $\quad \frac{d^{2n}}{dx^{2n}}f(x) = \frac{\Gamma(2n+2)}{\pi}\lim_{a\to 0}\frac{1}{a^{4n+1}}\int_{-\infty}^{\infty}\frac{P_{2n}(y;0,2n)}{(1+y^2)^{2n+1}}f(x+ay)dz$

Thus we arrive at a special case of (3.10) associated with the pseudo Jacobi polynomials.
When substituting $z = \frac{y}{ia}$ in the integral (7.7), the hypergeometric function can be written as an associated Legendre function of the second kind. This gives

(7.13) $\quad \frac{d^v}{dx^v}f(x) = \frac{e^{-(v+1/2)\pi i}a}{\pi}\int_{-\infty}^{\infty}\frac{Q_0^{v+1}(z)}{\left(\sqrt{z^2-1}\right)^{v+1}}f(x+iaz)dz$

**Example 2**.
Consider the function

(7.14) $\quad H^v(\omega) = (i\omega+0)^v e^{-\omega^2/2}$

The inverse Fourier transform is known (Erdélyi (1954) I. 121(23)) (note the factor $\sqrt{2\pi}$):

(7.15) $\quad h^v(y) = e^{-y^2/4}D_v(y) \qquad \mathrm{Re}\, v > -1$

Here $D_v(y)$ is a parabolic cylinder function. Although this is an old notation we use it here, mainly because of the relationship with the Hermite polynomials. With the known method one can derive for the fractional derivative

(7.16) $\quad \frac{d^v}{dx^v}f(x) = \lim_{\delta\to 0}\frac{1}{\delta^v}\int_{-\infty}^{\infty}e^{-y^2/4}D_v(y)f(x-\delta y)dy$

Formula (7.16) can also be computed when (3.10) will be applied to the Hermite polynomials. There follows

(7.17) $\quad W_\delta^{v,n}[f](x) = \frac{(-1)^n}{\sqrt{2\pi}\,\Gamma(n-v)}\frac{1}{\delta^v}\int_{-\infty}^{\infty}f(x+\delta y)\left[\int_{-\infty}^{y}(y-u)^{n-v-1}H_n(u)e^{-u^2}du\right]dy$

For the integral between the square brackets we use the Rodrigues formula. There follows

(7.18) $\quad I = \int_{-\infty}^{y}(y-u)^{n-v-1}H_n(u)e^{-u^2}du = (-1)^n\int_{-\infty}^{y}(y-u)^{n-v-1}\frac{d^n}{du^n}[e^{-u^2}]du$

Repeated partial integration $n$ times gives

(7.19) $\quad I = (-1)^n\frac{\Gamma(n-v)}{\Gamma(-v)}\int_{-\infty}^{y}(y-u)^{-v-1}e^{-u^2}du = (-1)^n\frac{\Gamma(n-v)}{\Gamma(-v)}e^{-y^2}\int_0^{\infty}u^{-1-v}e^{-u^2+2yu}du$

This integral is well-known (Erdélyi (1954) I. 313(13)). There is

(7.20) $$\int_0^\infty u^{-\nu-1} e^{-u^2+2yu} du = 2^{\nu/2} \Gamma(-\nu) e^{y^2/2} D_\nu\left(-y\sqrt{2}\right)$$

Substitution gives

(7.21) $$W_\delta^{\nu,n}[f](x) = \frac{1}{\sqrt{2\pi}} \frac{2^{\nu/2}}{\delta^\nu} \int_{-\infty}^\infty e^{-y^2/2} D_\nu\left(-y\sqrt{2}\right) f(x+\delta y) dy$$

Replacing $y$ by $-\dfrac{y}{\sqrt{2}}$ and $\delta$ by $\delta\sqrt{2}$ gives

(7.22) $$W_\delta^{\nu,n}[f](x) = \frac{1}{\sqrt{2\pi}} \frac{1}{\delta^\nu} \int_{-\infty}^\infty e^{-y^2/4} D_\nu(y) f(x-\delta y) dy$$

This result is the same as (7.16) when taking the limit.
When $\nu$ is an integer $n$ we obtain

(7.23) $$D_n(y) = 2^{-n/2} e^{-y^2/4} H_n\left(\frac{y}{\sqrt{2}}\right)$$

Substitution in (7.22) gives

(7.24) $$\frac{d^n}{dx^n} f(x) = \frac{2^{-n/2}}{\sqrt{2\pi}} \lim_{\delta \to 0} \frac{1}{\delta^n} \int_{-\infty}^\infty e^{-y^2/2} H_n\left(\frac{y}{\sqrt{2}}\right) f(x-\delta y) dy$$

Replacing $y$ by $y\sqrt{2}$ and $\delta$ by $\dfrac{\delta}{\sqrt{2}}$ there arises the formula for the Hermite derivative

(7.25) $$\frac{d^n}{dx^n} f(x) = \frac{1}{\sqrt{\pi}} \lim_{\delta \to 0} \frac{1}{\delta^n} \int_{-\infty}^\infty e^{-y^2/2} H_n(y) f(x-\delta y) dy$$

## 8 The fractional orthogonal difference.

In this section we treat the fractional difference. If we want to compute the fractional derivative of a function given at an equidistant set of points we can use instead of the fractional derivative the fractional difference. Then we have to work in the formula for the orthogonal derivative with discrete orthogonal polynomials.
As an analogue of the Riemann-Liouville fractional integral (2.1) we can define a fractional summation (see for example Kuttner (1957) or Isaacs (1963)) as

(8.1) $$I_\delta^{-\mu}[f](x) := \delta^\mu \sum_{k=0}^\infty \frac{(\mu)_k}{k!} f(x-k\delta) \qquad \delta > 0$$

In contrast with the Riemann-Liouville fractional integral, this summation has no a priori singularity for $\mu < 0$. If $f$ is a causal function then the upper bound for the summation should be $[x/\delta]$ and the sum is certainly well-defined. For $\mu \neq 0, -1, -2, \ldots$ we can write

(8.2) $$\frac{(\mu)_k}{k!} = \frac{1}{\Gamma(\mu)} \frac{\Gamma(\mu+k)}{\Gamma(1+k)} \sim \frac{k^{\mu-1}}{\Gamma(\mu)}$$

as $k \to \infty$. Hence the infinite sum converges if $f(x) = O(|x|^{-\lambda})$ as $x \to \infty$ with $\lambda > \operatorname{Re}\mu$. For $\mu = -n$ ($n$ a nonnegative integer) we get

(8.3a) $$I_\delta^n[f](x) = \delta^{-n} \sum_{k=0}^\infty (-1)^k \binom{n}{k} f(x-k\delta) = (D_\delta^-)^n[f](x)$$

where

(8.3b) $$D_\delta^-[f](x) := \frac{f(x) - f(x-\delta)}{\delta}$$

Note that

(8.4) $$D_\delta^-\left[I_\delta^{-\mu}[f]\right](x) = I_\delta^{-\mu+1}[f](x)$$

More generally we have

(8.5) $\quad I_\delta^{-\mu}\left[I_\delta^{-\nu}[f]\right](x) = I_\delta^{-(\mu+\nu)}[f](x)$

For $\mu > 0$ the fractional summation $I_\delta^{-\mu}[f]$ formally approximates the Riemann-Liouville integral $R^{-\mu}[f]$. Indeed, putting $y = k\delta$ we can rewrite the formula for $I_\delta^{-\mu}[f](x)$ as

(8.6) $\quad I_\delta^{-\mu}[f](x) = \dfrac{\delta^\mu}{\Gamma(\mu)} \sum\limits_{y \in \delta \mathbb{Z}_{\geq 0}} \dfrac{\Gamma\left(\mu + \frac{y}{\delta}\right)}{\Gamma\left(1 + \frac{y}{\delta}\right)} f(x-y) = \dfrac{\delta}{\Gamma(\mu)} \sum\limits_{y \in \delta \mathbb{Z}_{\geq 0}} \dfrac{y^{\mu-1}\Gamma\left(\mu + \frac{y}{\delta}\right)}{\left(\frac{y}{\delta}\right)^{\mu-1}\Gamma\left(1 + \frac{y}{\delta}\right)} f(x-y)$

Formula (5.11.12) of Olver et al. (2010) gives if $\delta \downarrow 0$:

(8.7) $\quad \lim\limits_{\delta \downarrow 0} \dfrac{\Gamma\left(\frac{y}{\delta} + \mu\right)}{\Gamma\left(\frac{y}{\delta} + 1\right)\left(\frac{y}{\delta}\right)^{\mu-1}} = 1$

Formally we have:

(8.8) $\quad \lim\limits_{\delta \downarrow 0} \delta \sum\limits_{y \in \delta \mathbb{Z}_{\geq 0}} g_\delta(y) = \int_0^\infty g(y) dy$

with $g_\delta(y) \to g(y)$ as $\delta \downarrow 0$. Application of (8.7) and (8.8) to (8.6) yields:

(8.9) $\quad \lim\limits_{\delta \downarrow 0} I_\delta^{-\mu}[f](x) = \dfrac{1}{\Gamma(\mu)} \int_0^\infty y^{\mu-1} f(x-y) dy$

This is the fractional integral of Riemann-Liouville. So it is shown that the fractional difference formally tends to the fractional integral in the limit as $\delta \downarrow 0$.

Remark 8.1. Grünwald (1867) and Letnikov (1868) developed a theory of discrete fractional calculus for causal functions ($f(x) = 0$ for $x < 0$). For a good description of their theory see (see Oldham, K.B., Spanier, J. 1974). They started with the backwards difference

(8.10) $\quad \Delta^-_\delta[f](x) = \dfrac{f(x) - f(x - \delta)}{\delta}$

and next they define the fractional derivative as

(8.11) $\quad \dfrac{d^\nu}{dx^\nu} f(x) := \dfrac{1}{\Gamma(-\nu)} \lim\limits_{\delta \to 0, N \to \infty} \dfrac{1}{\delta^\nu} \sum\limits_{m=0}^{N-1} \dfrac{\Gamma(m-\nu)}{\Gamma(m+1)} f(x - k\delta)$

where $N = \left[\frac{x}{\delta}\right]$. This definition has the advantage that $\nu$ is fully arbitrary. No use is made of derivatives or integrals. The disadvantage is that computation of the limit is often very difficult. However it can be shown that (8.11) is the same as the Riemann-Liouville derivative for causal functions for $\operatorname{Re}\nu > 0$ and $f$ is $n$ times differentiable with $n > \operatorname{Re}\nu$. At least formally this follows also by our approach if we write $I_\delta^\nu[f](x) = (D_\delta^-)^n[I_\delta^{\nu-n}[f]](x)$ and then use (8.9) and the fact that (8.3b) approaches the derivative as $\delta \downarrow 0$.

Comparing (8.11) with (8.1) the differences are the upper bound and the limit. For practical calculations formula (8.11) is often used without the limit. In the next section we will develop the fractional orthogonal difference for the Hahn polynomials. In section 10 we will see that the definition of Grünwald-Letnikov for the fractional difference is a special case of the definition of the fractional orthogonal difference for the Hahn polynomials for causal functions.

For obtaining the formula for the fractional orthogonal difference we let $\{p_n\}_{n=0}^N$ be a system of orthogonal polynomials with respect to weights $w(j)$ on points $j$ ($j = 0, 1, \ldots, N$). Then the approximate orthogonal derivative (3.1b) takes the form

(8.12) $\quad D_\delta^n[f](x) = \dfrac{k_n n!}{h_n} \dfrac{1}{\delta^n} \sum\limits_{j=0}^N f(x + j\delta) w(j) p_n(j) \qquad n = 0, 1, \ldots, N$

Let $\operatorname{Re}\nu < n$ with $n$ a positive integer. The fractional derivative $R^\nu[f](x) = D^n\left[R^{\nu-n}[f]\right](x)$ can be formally

approximated by

(8.13) $\quad I_{\delta,n,w}^{\nu}[f](x) := D_{\delta}^{n}\left[I_{\delta}^{\nu-n}[f]\right](x)$

as $\delta \downarrow 0$. This is a motivation to compute $I_{\delta,n,w}^{\nu}[f]$ more explicitly, in particular for special choices of the weights $w$. Substitution of (8.12) and (8.1) in (8.13) gives:

(8.14) $\quad I_{\delta,n,w}^{\nu}[f](x) = \dfrac{k_n n!}{h_n} \dfrac{1}{\delta^{\nu}} \sum_{k=0}^{\infty} \sum_{j=0}^{N} \dfrac{(n-\nu)_k}{k!} f(x+(j-k)\delta) w(j) p_n(j)$

For the double sum we can write

(8.15) $\quad \sum_{k=0}^{\infty}\sum_{j=0}^{N} g(k,j) = \sum_{j=0}^{N}\sum_{k=0}^{\infty} g(k,j) = \sum_{j=0}^{N}\sum_{m=-j}^{\infty} g(j+m,j) = \sum_{j=0}^{N}\sum_{m=1}^{\infty} g(j+m,j) + \sum_{j=0}^{N}\sum_{m=0}^{j} g(j-m,j) =$

$\quad = \sum_{m=1}^{\infty}\sum_{j=0}^{N} g(j+m,j) + \sum_{m=0}^{N}\sum_{j=m}^{N} g(j-m,j)$

Application of this formula to (8.14) yields

(8.16) $\quad I_{\delta,n,w}^{\nu}[f] = \dfrac{k_n}{h_n}\Gamma(n+1)\dfrac{1}{\delta^{\nu}}\sum_{m=1}^{\infty} f(x-m\delta)\left[\sum_{j=0}^{N}\dfrac{(n-\nu)_{j+m}}{(j+m)!} p_n(j)w(j)\right] +$

$\quad + \dfrac{k_n}{h_n}\Gamma(n+1)\dfrac{1}{\delta^{\nu}}\sum_{m=0}^{N} f(x+m\delta)\left[\sum_{j=m}^{N}\dfrac{(n-\nu)_{j-m}}{(j-m)!} p_n(j)w(j)\right]$

The summations inside the square brackets are in principle known after a choice of the discrete orthogonal polynomials $p_n(x)$. Then this formula can be used for approximating the fractional orthogonal derivative.

Remark 8.2. When taking the Fourier transform of the fractional summation one can show that for $\delta \downarrow 0$ this Fourier transform formally tends to the Fourier transform of the fractional integral. To this purpose we write (8.1) as

(8.17) $\quad I_{\delta}^{-\mu}[f](n\delta) = \delta^{\mu}\sum_{k=0}^{\infty}\dfrac{(\mu)_k}{k!} f\big((n-k)\delta\big)$

Taking the discrete Fourier transform gives

(8.18) $\quad \sum_{n=-\infty}^{\infty} I_{\delta}^{-\mu}[f](n\delta)e^{in\delta\omega} = \delta^{\mu}\sum_{n=-\infty}^{\infty} e^{in\delta\omega}\sum_{k=0}^{\infty}\dfrac{(\mu)_k}{k!} f\big((n-k)\delta\big)$

while working formally, we interchange the summations on the right. This gives:

(8.19) $\quad \sum_{n=-\infty}^{\infty} I_{\delta}^{-\mu}[f](n\delta)e^{in\delta\omega} = \delta^{\mu}\sum_{k=0}^{\infty}\dfrac{(\mu)_k}{k!} e^{ik\delta\omega}\sum_{n=-\infty}^{\infty} f\big((n-k)\delta\big)e^{i(n-k)\delta\omega} = \delta^{\mu}\sum_{k=0}^{\infty}\dfrac{(\mu)_k}{k!} e^{ik\delta\omega}\sum_{n=-\infty}^{\infty} f(n\delta)e^{in\delta\omega}$

The first summation is well known. We obtain

(8.20) $\quad \sum_{n=-\infty}^{\infty} I_{\delta}^{-\mu}[f](n\delta)e^{in\delta\omega} = \left(\dfrac{1-e^{i\delta\omega}}{\delta}\right)^{-\mu}\sum_{n=-\infty}^{\infty} f(n\delta)e^{in\delta\omega} = \omega^{-\mu}\left(\dfrac{\sin(\delta\omega/2)}{\delta\omega/2}\right)^{-\mu} e^{-i\mu(\delta\omega-\pi)/2}\sum_{n=-\infty}^{\infty} f(n\delta)e^{in\delta\omega}$

Multiply both sides with $\delta$ and take the formal limit for $\delta \downarrow 0$. We obtain the well known formula

(8.21) $\quad \int_{-\infty}^{\infty} R^{-\mu}[f](x)e^{ix\omega}dx = e^{i\pi\mu/2}\omega^{-\mu}\int_{-\infty}^{\infty} f(x)e^{ix\omega}dx$

## 9 The fractional orthogonal difference for the discrete Hahn polynomials.

In this section we substitute the discrete Hahn polynomials $Q_n(x;\alpha,\beta,N)$ in formula (8.14) to calculate the fractional difference. We do this analogous to Jacobi polynomials. For the Hahn polynomials we use the definition in Koekoek et al. (2010). They gave

(9.1) $\quad p_n(j) = Q_n(j;\alpha,\beta,N) := {}_3F_2\left(\begin{array}{c} -j,-n,n+\alpha+\beta+1 \\ \alpha+1,-N \end{array};1\right)$

For the Hahn polynomials there are the following properties:

(9.2a) $\quad n = 0,1,2,\ldots,N$

(9.2b) $\quad w(j;\alpha,\beta,N) = \dfrac{(\alpha+1)_j(\beta+1)_{N-j}}{j!(N-j)!}$

Substitution in (8.14) gives:

(9.3a) $\quad I^\nu_{\delta,n,w}[f] = H_1(n)\left[\sum_{m=1}^{\infty} f(x-m\delta)J_1(m) + \sum_{m=0}^{N} f(x+m\delta)J_2(m)\right]$

with

(9.3b) $\quad J_1(m) = \sum_{j=0}^{N} \dfrac{(n-\nu)_{j+m}}{(j+m)!}(wQ_n)(j;\alpha,\beta,N)$

(9.3c) $\quad J_2(m) = \sum_{j=m}^{N} \dfrac{(n-\nu)_{j-m}}{(j-m)!}(wQ_n)(j;\alpha,\beta,N) = \sum_{j=0}^{N-m} \dfrac{(n-\nu)_j}{j!}(wQ_n)(j+m;\alpha,\beta,N)$

(9.3d) $\quad H_1(n) = (-1)^n \dfrac{\Gamma(2n+\alpha+\beta+2)\Gamma(\beta+1)\Gamma(n+1)}{\Gamma(N+n+\alpha+\beta+2)\Gamma(n+\beta+1)} \dfrac{1}{\delta^\nu}$

To compute the sums in (9.3b) and (9.3c) we use partial summation. The general formulas are (5.17) and (5.19). For applying these formula for the summation in (9.3b) we use (5.14). Then

(9.4) $\quad J_1(m) = \dfrac{(-1)^n(\beta+1)_n}{(-N)_n}\sum_{j=0}^{N} \dfrac{(n-\nu)_{j+m}}{(j+m)!}\nabla^n_j\left(\dfrac{(\alpha+n+1)_j}{j!}\dfrac{(\beta+n+1)_{N-n-j}}{(N-n-j)!}\right)$

Suppose

(9.5a) $\quad f(j) := \dfrac{(\alpha+n+1)_j}{j!}\dfrac{(\beta+n+1)_{N-n-j}}{(N-n-j)!}$

(9.5b) $\quad g(j) := \dfrac{(n-\nu)_{j+m}}{(j+m)!} = \dfrac{(j+m+1)_{n-\nu-1}}{\Gamma(n-\nu)}$

and seeing that $f(N-n) = f(-1) = 0$ we use (5.18) to get

(9.6) $\quad J_1(m) = \dfrac{(\beta+1)_n}{\Gamma(n-\nu)(-N)_n}\sum_{j=0}^{N-n}\dfrac{(\alpha+n+1)_j}{j!}\dfrac{(\beta+n+1)_{N-n-j}}{(N-n-j)!}\nabla^n_j\Big((j+m+n+1)_{n-\nu-1}\Big)$

For the second fraction inside the summation we use (Slater (1966), I.30)

(9.7) $\quad \Gamma(a-j) = \dfrac{\Gamma(a)(-1)^j}{(1-a)_j}$

Then there follows

(9.8) $\quad \dfrac{(\beta+n+1)_{N-n-j}}{(N-n-j)!} = \dfrac{(\beta+n+1)_{N-n}}{(N-n)!}\dfrac{(n-N)_j}{(-\beta-N)_j}$

For the third fraction inside the summation we use the formula

(9.9) $\quad \nabla^n\Big((j+a)_b\Big) = (a)_{b-n}(b-n+1)_n\dfrac{(a+b-n)_j}{(a)_j}$

Then we obtain

(9.10) $\quad \nabla_j^n\left((j+m+n+1)_{n-\nu-1}\right) = (-\nu)_n \dfrac{\Gamma(m+n-\nu)}{\Gamma(m+n+1)} \dfrac{(m+n-\nu)_j}{(m+n+1)_j}$

After substitution in $J_1(m)$ we get

(9.11) $\quad J_1(m) = \dfrac{1}{\Gamma(-\nu)(-N)_n} \dfrac{(1+\beta)_N}{(N-n)!} \dfrac{\Gamma(m+n-\nu)}{\Gamma(m+n+1)} \sum_{j=0}^{N-n} \dfrac{(\alpha+n+1)_j(m+n-\nu)_j(n-N)_j}{(m+n+1)_j(-\beta-N)_j} \dfrac{1}{j!}$

Then $J_1(m)$ can be written as a hypergeometric function

(9.12) $\quad J_1(m) = \dfrac{1}{(-N)_n} \dfrac{(1+\beta)_N}{(N-n)!} \dfrac{(-\nu)_{m+n}}{(m+n)!} \; _3F_2\!\left(\begin{array}{c} n-N, n+\alpha+1, m+n-\nu \\ -N-\beta, m+n+1 \end{array}; 1\right)$

Observing the hypergeometric function we will see that this function represents a series with a finite number of terms (if $\beta$ is not an integer with $-N-\beta \neq 0$). With $\alpha = \beta = 0$ and $\operatorname{Re}(\alpha) < m$ there follows

(9.13) $\quad J_1(m) = \dfrac{(N+1-n)_n}{(-N)_n} \dfrac{(-\nu)_{m+n}}{(m+n)!} \; _3F_2\!\left(\begin{array}{c} n-N, n+1, m+n-\nu \\ -N, m+n+1 \end{array}; 1\right)$

Observing the hypergeometric function we will see that this function represents an infinite series.

To compute $J_2(m)$ we write formula (9.3c) with (5.14) as

(9.14) $\quad J_2(m) = \dfrac{(-1)^n (\beta+1)_n}{(-N)_n} \sum_{j=0}^{N-m} \dfrac{(n-\nu)_j}{j!} \nabla_j^n\!\left( \dfrac{(\alpha+n+1)_{j+m}}{(j+m)!} \dfrac{(\beta+n+1)_{N-n-m-j}}{(N-n-m-j)!} \right)$

Suppose

(9.15a) $\quad f(j) = \dfrac{(\alpha+n+1)_{j+m}}{(j+m)!} \dfrac{(\beta+n+1)_{N-n-m-j}}{(N-n-m-j)!}$

(9.15b) $\quad g(j) = \dfrac{(n-\nu)_j}{j!}$

Then we rewrite the formula for the partial summation (5.17) as

(9.16) $\quad \sum_{j=-1}^{N-m} (\nabla f)(j) g(j) = f(N-m) g(N-m) - f(-2) g(-1) - \sum_{j=-1}^{N-m-1} f(j) \nabla_j^n(j+1)$

In particular, if $f(N-m) = g(-1) = 0$ then

(9.17) $\quad \sum_{j=0}^{N-m} (\nabla f)(j) g(j) = - \sum_{j=-1}^{N-m-1} f(j)(\nabla g)(j+1) = - \sum_{j=0}^{N-m} f(j-1)(\nabla g)(j)$

Applying this formula $n$ times gives

(9.18) $\quad \sum_{j=0}^{N-m} (\nabla f)^n(j) g(j) = (-1)^n \sum_{j=0}^{N-m} f(j-n)(\nabla g)^n(j)$

Substitution in (9.14) gives

(9.19) $\quad J_2(m) = \dfrac{(\beta+1)_n}{\Gamma(n-\nu)(-N)_n} \sum_{j=0}^{N-m} \dfrac{(\alpha+n+1)_{j+m-n}}{(j-n+m)!} \dfrac{(\beta+n+1)_{N-m-j}}{(N-m-j)!} \nabla^n\!\left((j+1)_{n-\nu-1}\right)$

For the $\nabla$ function we use formula (9.9). We obtain

(9.20) $$\nabla^n\left((j+1)_{n-\nu-1}\right) = \Gamma(n-\nu)\frac{(-\nu)_j}{j!}$$

After substitution in $J_2(m)$ we get

(9.21) $$J_2(m) = \frac{(\beta+1)_n}{\Gamma(n-\nu)(-N)_n}\sum_{j=0}^{N-m}\frac{(\alpha+n+1)_{j+m-n}}{(j-n+m)!}\frac{(\beta+n+1)_{N-m-j}}{(N-m-j)!}\frac{(-\nu)_j}{j!}$$

For $j = 0$ the factor $(m - n + 1)$ could be less or equal then zero. Then the terms concerning this condition should not be taken with. In this case the summation should be reversed. We set $j = N - m - i$. This gives

(9.22) $$J_2(m) = \frac{(\beta+1)_n}{\Gamma(n-\nu)(-N)_n}\sum_{i=0}^{\min(N-m,N-n)}\frac{(\alpha+n+1)_{N-n-i}}{(N-n-i)!}\frac{(\beta+n+1)_i}{i!}\frac{(-\nu)_{N-m-i}}{(N-m-i)!}$$

After some manipulations with the Pochhammer symbols, $J_2(m)$ can be written as a hypergeometric function.

(9.23) $$J_2(m) = \frac{(\beta+1)_n}{(-N)_n}\frac{(\alpha+n+1)_{N-n}}{(N-n)!}\frac{(-\nu)_{N-m}}{(N-m)!}\;_3F_2\left(\begin{array}{c}n-N,\beta+n+1,-N+m\\-\alpha-N,\nu-N+m+1\end{array};1\right)$$

With $\alpha = \beta = 0$ there follows after again some manipulations of the Pochhammer symbols

(9.24) $$J_2(m) = \frac{(-1)^n}{\Gamma(-\nu)}\frac{\Gamma(N-m-\nu)}{(N-m)!}\;_3F_2\left(\begin{array}{c}n-N,n+1,-N+m\\-N,\nu-N+m+1\end{array};1\right)$$

It can be shown that the hypergeometric function represents a finite series. Then with (Olver et al (2010),16.4.11) we can rewrite the hypergeometric function. This results in

(9.25) $$J_2(m) = \frac{\Gamma(\nu-2n)}{\Gamma(\nu+1-n)\Gamma(-\nu)}\frac{\Gamma(N-m-\nu+n+1)}{(N-m)!}\;_3F_2\left(\begin{array}{c}-n,n+1,-m\\-N,\nu+1-n\end{array};1\right)$$

For the approximate fractional difference there results from (9.3), (9.12) and (9.23)

(9.26) $$I_{\delta,n}^\nu[f] = \frac{\Gamma(n+1)}{\Gamma(N+1)}\frac{\Gamma(2n+\alpha+\beta+2)}{\Gamma(N+n+\alpha+\beta+2)}\frac{\Gamma(N+1+\beta)}{\Gamma(n+1+\beta)}\frac{1}{\delta^\nu}$$
$$\sum_{m=1}^{\infty}f(x-m\delta)\frac{(-\nu)_{m+n}}{(m+n)!}\;_3F_2\left(\begin{array}{c}n-N,n+\alpha+1,m+n-\nu\\-N-\beta,m+n+1\end{array};1\right) +$$
$$+\frac{\Gamma(n+1)}{\Gamma(N+1)}\frac{\Gamma(2n+\alpha+\beta+2)}{\Gamma(N+n+\alpha+\beta+2)}\frac{\Gamma(\alpha+N+1)}{\Gamma(\alpha+n+1)}\frac{1}{\delta^\nu}$$
$$\sum_{m=0}^{N}f(x+m\delta)\frac{(-\nu)_{N-m}}{(N-m)!}\;_3F_2\left(\begin{array}{c}n-N,\beta+n+1,-N+m\\-\alpha-N,\nu-N+m+1\end{array};1\right)$$

For $\alpha = \beta = 0$ there remains

(9.27) $$I_{\delta,n}^\nu[f] = \frac{\Gamma(2n+2)}{\Gamma(N+n+2)}\frac{1}{\Gamma(-\nu)}\frac{1}{\delta^\nu}\sum_{m=1}^{\infty}f(x-m\delta)\frac{\Gamma(m+n-\nu)}{(m+n)!}\;_3F_2\left(\begin{array}{c}n-N,n+1,m+n-\nu\\-N,m+n+1\end{array};1\right) +$$
$$+(-1)^n\frac{\Gamma(2n+2)}{\Gamma(N+n+2)}\frac{1}{\Gamma(-\nu)}\frac{\Gamma(\nu-2n)}{\Gamma(\nu+1-n)}\frac{1}{\delta^\nu}$$
$$\sum_{m=0}^{N}f(x+m\delta)\frac{\Gamma(N-m-\nu+n+1)}{(N-m)!}\;_3F_2\left(\begin{array}{c}-n,n+1,-m\\-N,\nu+1-n\end{array};1\right)$$

For $n = 1$ there follows

(9.28) $$I_{\delta,1}^\nu[f] = \frac{6}{\Gamma(N+3)}\frac{1}{\Gamma(-\nu)}\frac{1}{\delta^\nu}\sum_{m=1}^{\infty}f(x-m\delta)\frac{\Gamma(m+1-\nu)}{\Gamma(m+2)}\;_3F_2\left(\begin{array}{c}1-N,2,m+1-\nu\\-N,m+2\end{array};1\right) -$$
$$-\frac{6}{\Gamma(N+3)}\frac{1}{\Gamma(-\nu)}\frac{\Gamma(\nu-2)}{\Gamma(\nu)}\frac{1}{\delta^\nu}\sum_{m=0}^{N}f(x+m\delta)\frac{\Gamma(N-m-\nu+2)}{\Gamma(N-m+1)}\;_3F_2\left(\begin{array}{c}-1,2,-m\\-N,\nu\end{array};1\right)$$

The first hypergeometric function can be written as

(9.29) $$\,_3F_2\left(\begin{array}{c}1-N,2,m+1-\nu\\-N,m+2\end{array};1\right)=\frac{\Gamma(m+2)}{(-N)\Gamma(m+1-\nu)}\sum_{k=0}^{N-1}(k-N)(k+1)\frac{\Gamma(k+m+1-\nu)}{\Gamma(k+m+2)}$$

Making use of

(9.30) $$(k-N)(k+1)=(k+m+1)(k+m)-(N+2m)(k+m+1)+m(N+m+1)$$

and of (Erdélyi (1953) I. 2.5 (16))

(9.31) $$\sum_{k=0}^{N-1}\frac{\Gamma(k+a)}{\Gamma(k+b)}=\frac{1}{(a-b+1)}\left[\frac{\Gamma(N+a)}{\Gamma(N+b-1)}-\frac{\Gamma(a)}{\Gamma(b-1)}\right]\quad\left(a+1\neq b,\ N>0\right)$$

there remains

(9.32) $$\,_3F_2\left(\begin{array}{c}1-N,2,m+1-\nu\\-N,m+2\end{array};1\right)=$$
$$=\frac{\Gamma(m+2)\Gamma(-\nu)}{N\Gamma(m+1-\nu)\Gamma(3-\nu)}\left[(2N+2m-2\nu-N\nu+2)\frac{\Gamma(m-\nu+1)}{\Gamma(m)}-(2m+N\nu)\frac{\Gamma(N+m-\nu+2)}{\Gamma(N+m+1)}\right]$$

For the second hypergeometric functions there yields

(9.33) $$\,_3F_2\left(\begin{array}{c}-1,2,-m\\-N,\nu\end{array};1\right)=\frac{(N\nu-2m)}{N\nu}$$

After substitution in (9.28) there remains

(9.34) $$I_{\delta,1}^{\nu}[f](x)\left[\frac{6}{N\Gamma(N+3)\Gamma(3-\nu)}\frac{1}{\delta^{\nu}}\right]^{-1}=$$
$$=\sum_{m=1}^{\infty}f(x-m\delta)\left[(2N+2m-2\nu-N\nu+2)\frac{\Gamma(m-\nu+1)}{\Gamma(m)}-(2m+N\nu)\frac{\Gamma(N+m-\nu+2)}{\Gamma(N+m+1)}\right]+$$
$$+\sum_{m=0}^{N}f(x+m\delta)\frac{\Gamma(N-m-\nu+2)}{\Gamma(N-m+1)}(2m-N\nu)$$

## 10  The transfer function for the approximate fractional orthogonal difference.

Just as for the fractional orthogonal derivative it is possible to compute the transfer function for the fractional orthogonal difference. In that case we can apply Theorem 2.1. where we take another function for the function $g$. In this case we use (8.1)

(10.1) $$g=I_{\delta}^{\mu}[f](x)=\delta^{\mu}\sum_{k=0}^{\infty}\frac{(\mu)_k}{k!}f(x-k\delta)\qquad \delta>0$$

Fourier transform of this function gives

(10.2) $$G(\omega)=\frac{1}{\sqrt{2\pi}}\int_{-\infty}^{\infty}e^{i\omega x}g(x)dx=\delta^{\mu}\sum_{k=0}^{\infty}\frac{(\mu)_k}{k!}\frac{1}{\sqrt{2\pi}}\int_{-\infty}^{\infty}e^{i\omega x}f(x-k\delta)dx=\delta^{\mu}\sum_{k=0}^{\infty}\frac{(\mu)_k}{k!}e^{i\omega k\delta}F(\omega)$$

The summation is known and there follows

(10.3) $$G(\omega)=\left(\frac{1-e^{i\omega\delta}}{\delta}\right)^{-\mu}F(\omega)$$

From (8.10), (8.11), (5.1) and (10.3) we obtain that if $g=I_{\delta,n,w}^{\nu}[f]$, then

(10.4) $$H(\omega)=\frac{G(\omega)}{F(\omega)}=\left(\frac{1-e^{i\omega\delta}}{\delta}\right)^{\nu-n}\frac{k_n n!}{h_n\delta^n}\sum_{x=0}^{N}p_n(x)w(x)e^{-ix\delta\omega}$$

For the Hahn polynomials this becomes

(10.5) $$H(\omega)=\left(\frac{1-e^{i\omega\delta}}{\delta}\right)^{\nu-n}\frac{k_n n!}{h_n\delta^n}S$$

The summation $S$ is defined by (5.17) and computed in (5.25). We obtain

$$(10.6) \qquad H(\omega) = \left(\frac{1-e^{i\omega\delta}}{\delta}\right)^{\nu} \frac{k_n n!}{h_n} \frac{(\beta+1)_N}{N!} F\left(\begin{array}{c} -N+n, \alpha+n+1 \\ -\beta-N \end{array}; e^{-i\delta\omega}\right)$$

Substitution of (5.27) yields

$$(10.7) \qquad H(\omega) = (-1)^n \left(\frac{1-e^{i\omega\delta}}{\delta}\right)^{\nu} \frac{\Gamma(N+\beta+1)\Gamma(2n+\alpha+\beta+2)}{\Gamma(n+\beta+1)\Gamma(N+n+\alpha+\beta+2)} F\left(\begin{array}{c} -N+n, \alpha+n+1 \\ -\beta-N \end{array}; e^{-i\delta\omega}\right)$$

Comparing this formula with formula (5.28) we will see that the only difference between these formulas is the frequency dependent factor. In the next section we will discuss this fact.

With $\alpha = \beta = 0$ and $n = 1$ there remains

$$(10.8) \qquad H(\omega) = -\left(\frac{1-e^{i\omega\delta}}{\delta}\right)^{\nu} \frac{6\Gamma(N+1)}{\Gamma(N+3)} F\left(\begin{array}{c} 1-N, 2 \\ -N \end{array}; e^{-i\delta\omega}\right) =$$

$$= -\left(\frac{1-e^{i\omega\delta}}{\delta}\right)^{\nu} \frac{6\Gamma(N)}{\Gamma(N+3)} \sum_{k=0}^{N-1}(k-N)(k+1)e^{-i\delta k\omega} =$$

$$= \frac{(1-e^{i\omega\delta})^{\nu}}{\delta^{\nu}(e^{-i\delta\omega}-1)^3} \frac{6\Gamma(N)}{\Gamma(N+3)} \left(e^{-i\delta\omega(N+1)}(N-Ne^{-i\delta\omega}+2)-(2e^{-i\delta\omega}+Ne^{-i\delta\omega}-N)\right)$$

For $N = 1$ there remains

$$(10.9) \qquad H(\omega) = -\left(\frac{1-e^{i\omega\delta}}{\delta}\right)^{\nu}$$

Remark 10.1. Note that the transfer function in (10.7) is derived with a summation to infinite. For practical reasons one should take the summation to a finite number say $M$. Then we have to use the transfer function corresponding to (9.34). This is important because (9.34) is the formula to compute the approximate fractional Hahn derivative. Then the transfer function will differ from the ideal one. Taking the Fourier transform gives for the transfer function

$$(10.10) \qquad H(\omega)\left[\frac{6}{\Gamma(N+3)}\frac{1}{\delta^{\nu}}\right]^{-1} =$$

$$= \frac{1}{\Gamma(2-\nu)}e^{-i(N+1)\omega\delta}\sum_{m=1}^{M+N+1}e^{im\omega\delta}\frac{\Gamma(m-\nu+1)}{\Gamma(m)} + \frac{1}{\Gamma(2-\nu)}\sum_{m=1}^{M}e^{im\omega\delta}\frac{\Gamma(m-\nu+1)}{\Gamma(m)} -$$

$$- \frac{2}{N}\frac{1}{\Gamma(3-\nu)}e^{-iN\omega\delta}\sum_{m=1}^{M+N}e^{im\omega\delta}\frac{\Gamma(m-\nu+2)}{\Gamma(m)} + \frac{2}{N}\frac{1}{\Gamma(3-\nu)}\sum_{m=1}^{M}e^{im\omega\delta}\frac{\Gamma(m-\nu+2)}{\Gamma(m)}$$

For $M \to \infty$ this function becomes the same as (10.8).

If $f(x)$ is causal then (10.7) can be derived from a finite sum and we do not cut of the summation. Then we can use the transfer function (10.7) instead of (10.10).

Remark 10.2. For the filter for the fractional differentiation of Grünwald and Letnikov the approximate transfer function (with $N$ finite) can be simply derived from (8.21) by taking the Fourier transform. We first take the Fourier transform of (8.21) without the limits. This gives with $N = [x/\delta]$

$$(10.11) \qquad F[GL^{\nu}[f](x)](\omega) = \frac{1}{\Gamma(-\nu)}\int_0^{\infty}e^{i\omega x}\frac{1}{\delta^{\nu}}\sum_{k=0}^{[x/\delta]}\frac{\Gamma(k-\nu)}{\Gamma(k+1)}f(x-k\delta)dx$$

Because $f(x)$ is a causal function it follows that $x \geq k\delta$. Substitution gives

$$(10.12) \qquad H(\omega) = \frac{1}{\Gamma(-\nu)}\frac{1}{\delta^{\nu}}\sum_{k=0}^{\infty}\frac{\Gamma(k-\nu)}{\Gamma(k+1)}\int_{k\delta}^{\infty}e^{i\omega x}f(x-k\delta)dx =$$

$$= \frac{1}{\Gamma(-\nu)}\frac{1}{\delta^{\nu}}\sum_{k=0}^{\infty}\frac{\Gamma(k-\nu)}{\Gamma(k+1)}e^{i\omega k\delta}F(\omega) = \left(\frac{1-e^{i\omega\delta}}{\delta}\right)^{\nu}F(\omega)$$

For $\delta \to 0$ the transfer function becomes

$$(10.13) \quad H(\omega) = \lim_{\delta \to 0} \left( \frac{1 - e^{i\omega\delta}}{\delta} \right)^\nu = (-i\omega)^\nu$$

Because of section 3.11 from Diekema and Koornwinder (2012) the Grünwald-Letnikov derivative follows by taking $w(j) = 1$ and $N = n$ in (8.12) or $\alpha = \beta = 0$ and $N = n$ in (9.1). So (10.12) becomes a special case of (10.7).

## 11  Application of the theory to a fractional differentiating filter.

This section treats the application of the fractional derivative in the linear filter theory. In that theory a filter is described by three properties, e.g. the input signal, the so-called impulse response function of the filter and the output signal. These properties are described in the time domain. The output signal is the convolution integral of the input signal and the impulse response function of the filter. The linearity means that the output signal depends linearly of the input signal.

One distinguishes between continuous and discrete filters, corresponding to continuity and discreteness, respectively, of the signals involved. In the discrete case the output signal can be computed with a discrete filter equation.

In the usage of filters there are two important items for consideration. The first item is the computation of the output signal given the input signal. The second item is a qualification of the working of the filter.

In our opinion this second task should be preferably done in the frequency domain, where one can see the spectra of the signal and eventually the noise. In this domain one obtains the Fourier transform of the output signal by multiplication of the Fourier transform of the input signal with the transfer function.

In this Section we will give graphs of the absolute value of the transfer function associated with various approximate fractional derivatives introduced in this paper. From these graphs we will get insight how these approximations will work in practice. We think this kind of analysis, in particular when using log-log plots, is preferable to the analysis of filters in the time domain. As an example of the latter tpye of analysis, actually for the fractional Jacobi derivative, see a paper of Tenreiro Machado (2009a).

For the second item the transfer function of the filter should be computed. We treat here the approximate fractional Jacobi derivative and the approximate fractional Hahn derivative. For practical computations these become the approximate fractional Legendre derivative and the approximate fractional Gram derivative.

For the first item we can use the formula (4.18) for the Jacobi derivative and formula (9.26) for the Hahn derivative. For the second item we had to go to the frequency domain.

In the frequency domain we suppose that the Fourier transforms of the input signal $x(t)$ and the output signal $y(t)$ are $X(\omega)$ and $Y(\omega)$ and does exist. The Fourier transformation of the impulse response function $h(t)$ of the filter is the transfer function $H(\omega)$ of the filter. For this transfer function there is the definition:

$$(11.1) \quad H(\omega) := \frac{Y(\omega)}{X(\omega)}$$

The function $H(\omega)$ is a complex function. For a graphic display of this function one uses the modulus and the phase of the transfer function. In our case we want to see the property of the differentiation and then the graph of the modulus of the transfer function is preferred.

For a $n$-th order differentiator one can show

$$(11.2) \quad H(\omega) = (-i\omega)^n$$

In the frequency domain the graph of the modulus of the transfer function of this filter is a straight line with a slope which depends on $n$. See figure 1.

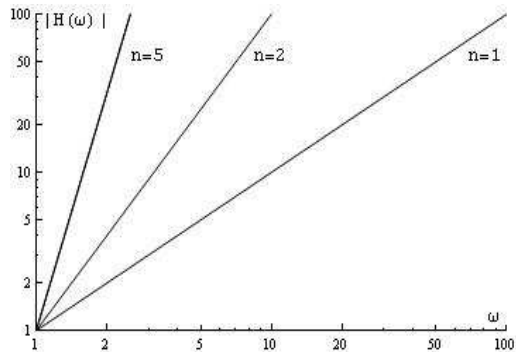

Figure 1: Moduli of the transfer function of a $n$-th order differentiator for $n = 1$, $n = 2$ and $n = 5$.

We see that for $\omega \to \infty$ the modulus of the transfer function $|H(\omega)|$ goes to $\infty$. This means that (because every system has some high frequency noise) this filter is instable. That is the reason why a differentiating filter should always have an integrating factor. This factor will be added to the filter so that for $\omega \to \infty$ the modulus of the transfer function goes to zero. The filter for the orthogonal derivatives has this property. These derivatives appears from an averaging process (least squares) and such a process will always give an integrating factor.

The formula for the transfer function of the approximate Legendre derivative is (Diekema and Koornwinder (2012) section 5):

$$(11.3) \quad H(\omega) = \frac{\Gamma(2n+2)}{2^n \Gamma(n+1)} \frac{1}{\delta^n} j_n(\omega\delta)$$

where the functions $j_n$ are the spherical Bessel functions (5.15). The graph of the modulus of the first order approximate Legendre derivative with $\delta = 1$ is given in figure 2. It is clear that for large $\omega$ the modulus of the transfer function goes to zero. For small $\omega$ the graph is a straight line with slope 1.

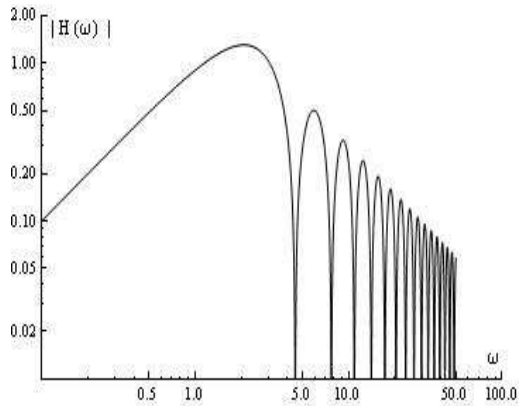

Figure 2: Modulus of the transfer function of the first order Legendre derivative.

In the case that the order of the differentiating filter is not an integer, the formula of the transfer function should be in analogue of (11.2)

$$(11.5) \quad H(\omega) = (-i\omega)^\nu$$

where $\nu$ can be an integer, a fractional or even a complex number. The graph of the modulus of the transfer function is given in figure 3.

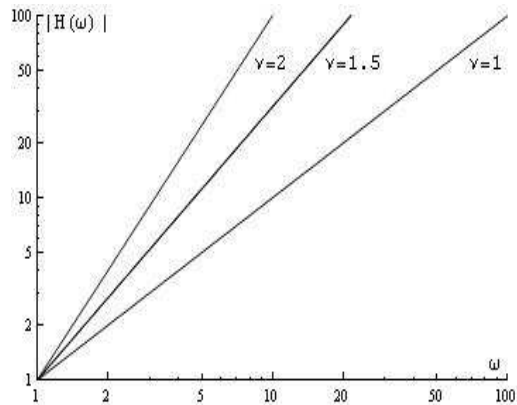

Figure 3: Moduli of the transfer functions of a fractional derivative with $\nu = 1$, $\nu = 1.5$ and $\nu = 2$.

Also in this case there is an instability. To prevent this instability one can use an approximate fractional Jacobi differentiating filter. For the approximate Jacobi derivative the transfer function is already computed in section 5. We repeat formula (5.5).

$$(11.6) \quad H(\omega) = (-i\omega)^\nu e^{-i\omega\delta} M(n+\alpha+1, 2n+\alpha+\beta+2; 2i\omega\delta)$$

With $\delta = 1$ the squared absolute value of the transfer function of the approximate fractional Jacobi derivative can be written as a series expansion as follows

(11.7) $$|H(\omega)|^2 = \sum_{m=0}^{\infty}(-4\omega^2)^m \sum_{k=0}^{2m} \frac{(n+\alpha+\beta+1)_k(n+\alpha+\beta+1)_{2m-k}}{(2n+\alpha+\beta+2)_k(2n+\alpha+\beta+1)_{2m-k}} \frac{(-1)^k}{k!(2m-k)!}$$

For small $\omega$ a first approximation of the modulus of the transfer function gives

(11.8) $$|H(\omega)| \sim \omega^\nu \sqrt{1 - \frac{4(n+\alpha+1)(n+\beta+1)}{(2n+\alpha+\beta+2)^2(2n+\alpha+\beta+3)}\omega^2}$$

It is clear that the formula is symmetric in $\alpha$ and $\beta$. The choices of $\alpha$ and $\beta$ detect the cut-off frequency. For simplicity we look at the frequency (and not at the exact cut-off frequency) where the transfer function has a first maximum. From formula (11.8) the maximum frequency is

(11.9) $$\omega_{max} = \frac{(2n+\alpha+\beta+2)}{2}\sqrt{\frac{\nu(2n+\alpha+\beta+3)}{(n+\alpha+1)(n+\beta+1)(1+\nu)}}$$

To simplify this formula set $\alpha = \beta$. Then there remains

(11.10) $$\omega_{max} = \sqrt{\frac{\nu(2n+2\alpha+3)}{(1+\nu)}}$$

The shape of the curves does not change in principle. So if one wants the cut-off frequency as high as possible, then $\alpha$ and $\beta$ should be chosen as high as possible. For the case $\alpha = \beta = 0$ the formula for the modulus of the transfer function of the approximate fractional Jacobi derivative simplifies to the approximate fractional Legendre derivative (5.14)

(11.11) $$H(\omega) = \omega^\nu \frac{\Gamma(2n+2)}{2^n \Gamma(n+1)} \frac{1}{(\omega\delta)^n} j_n(\omega\delta)$$

where the functions $j_n$ are the spherical Bessel functions with $n-1 \leq \nu \leq n$. For $n = 1$ and $\delta = 1$ there remains (see 5.15)

(11.12) $$H(\omega) = 3\omega^{\nu-3}(\sin\omega - \omega\cos\omega)$$

The graph of the modulus of this transfer function is given in figure 4 for different values of $\nu$.

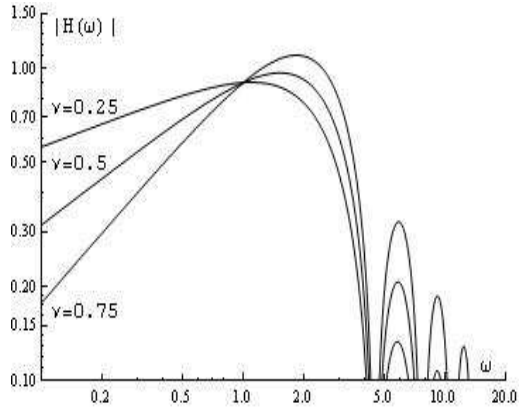

Figure 4: Moduli of the transfer functions of a fractional Legendre derivative.

Next we treat the discrete filter.

When we apply a discrete filter the input signal and the output signal are sampled with a sample frequency which is twice the maximum frequency of the input signal. So there is always a maximum frequency for the transfer function of the filter. What the transfer function will do above this maximum frequency is not important presupposing that for $\omega \to \infty$ the transfer function will go to zero. The input signal is known over $N$ points. So the input signal is an approximation of the true input signal.

The output signal can be computed for the $N$ points of the input signal. This is the reason that the output signal is always an approximation of the filtered true input signal. If the filter is a differentiator the output signal is an approximation of the derivative of the input signal. This approximation is very dependent of $N$.

For the fractional orthogonal derivative we use the transfer function of the approximate derivative of the discrete Hahn difference as derived in (10.7). For $\alpha = \beta = 0$ this transfer function becomes that of the fractional Gram derivative. In the following figure we take with $0 < v \leq 1$ and $n = 1$ (10.8). For $N = 1$ the modulus of the transfer function of the ideal filter (10.9) has a maximum at $\omega = \pi$. For different values of $N$ the modulus of the transfer function with $\delta = 1$ is given in the next figure.

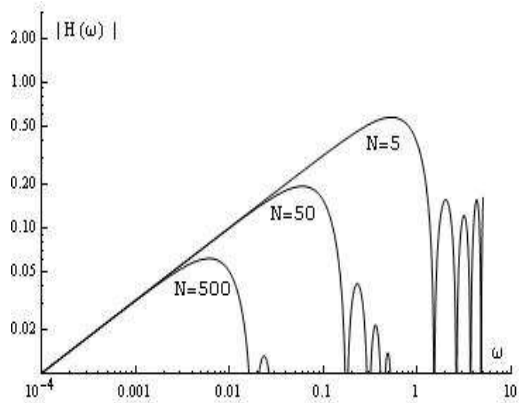

Figure 5: Moduli of the transfer function of the fractional Gram derivative with $v = 0.5$.

In practice (taking a finite number of points) we use the transfer function in (10.10). For $\delta = 1$, $v = 0.5$, $N = 7$ and different values of $M$ the modulus of this transfer function is shown in figure 6.

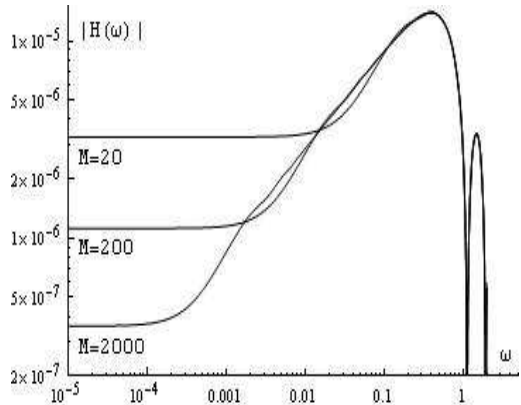

Figure 6. Moduli of the transfer functions of the fractional Hahn derivative with $N = 7$ and $v = 0.5$.

We see that for $\omega \to 0$ the absolute value of the transfer function tends to a constant value. Taking $\omega = 0$ in (10.10) there follows

$$(11.13) \quad H(0)\left[\frac{6}{\Gamma(N+3)}\frac{1}{\delta^v}\right]^{-1} = \frac{1}{\Gamma(2-v)}\sum_{m=1}^{M+N+1}\frac{\Gamma(m-v+1)}{\Gamma(m)} + \frac{1}{\Gamma(2-v)}\sum_{m=1}^{M}\frac{\Gamma(m-v+1)}{\Gamma(m)} -$$
$$- \frac{2}{N}\frac{1}{\Gamma(3-v)}\sum_{m=1}^{M+N}\frac{\Gamma(m-v+2)}{\Gamma(m)} + \frac{2}{N}\frac{1}{\Gamma(3-v)}\sum_{m=1}^{M}\frac{\Gamma(m-v+2)}{\Gamma(m)}$$

For the summations it is to prove

$$(11.14) \quad \sum_{m=1}^{K}\frac{\Gamma(m-a)}{\Gamma(m)} = \frac{\Gamma(K+1-a)}{\Gamma(K)(1-a)} \qquad a \neq 1$$

Then there remains

$$(11.15) \quad H(0)\left[\frac{6}{N\Gamma(N+3)\Gamma(4-v)}\frac{1}{\delta^v}\right]^{-1} =$$
$$= (N - 2M - Nv)\frac{\Gamma(M+N-v+3)}{\Gamma(M+N+1)} + \left((3-v)N + 2(M-v+2)\right)\frac{\Gamma(M-v+2)}{\Gamma(M)}$$

The lower bound of the frequency for which the filter does a fractional differentiation can be defined with the help of the following equation

$$\text{(11.16)} \quad \omega_l = \left( \frac{6(N - 2M - Nv)}{N\Gamma(N+3)\Gamma(4-v)} \frac{\Gamma(M+N-v+3)}{\Gamma(M+N+1)} + \frac{6\big((3-v)N + 2(M-v+2)\big)}{N\Gamma(N+3)\Gamma(4-v)} \frac{\Gamma(M-v+2)}{\Gamma(M)} \right)^{1/v}$$

In practice we should take this frequency 10 times higher then computed with (11.16).
For the maximum frequency we should compute the value of the frequency for the first maximum of the absolute value of (10.10). A tedious computation leads to

$$\text{(11.17)} \quad \omega_{\max} \approx 2\sqrt{6} \frac{\sqrt{(1-v)(6N + v + 6Nv + N^2 v + N^2 + 9)}}{\delta(6N + v + 6Nv + N^2 v + N^2 + 9)}$$

So we can define a bandwidth $B$ of the filter which is equal to $B = \omega_{\max} - \omega_l$.

For $M \to \infty$ we had to take the limits for the ratios of the gamma functions. Olver et al. (2010) formula 5.11.13 gives

$$\text{(11.18)} \quad \frac{1}{z^{a-b}} \frac{\Gamma(z+a)}{\Gamma(z+b)} = 1 + \frac{1}{2z}(a-b)(a+b-1) + $$
$$+ \frac{1}{24z^2}(a-b-1)(a-b)(3(a+b-1)^2 - (a-b+1)) + O(z^{-3})$$

Substitution in (11.15) gives with $\delta = 1$

$$\text{(11.19)} \quad H(0) = \frac{1}{M^v} \frac{(12N - 10v - 18Nv + 11v^2 - 3v^3 + 6Nv^2 - 6N^2 v + 6N^2)}{2N(N+2)\Gamma(N+1)(3-v)\Gamma(1-v)} \left(1 + O\!\left(\frac{1}{M}\right)\right)$$

Because $v > 0$ this constant goes to zero and the transfer function approximates the ideal case for low frequencies.
With this transfer function it is seen that the filter works only well when the frequencies of the input signal lie inside the band-pass of the filter. The maximum frequency of this band-pass is mainly dependent on the number $N$ and the minimum frequency on the number $M$.
Many authors (e.g. Diethelm et al. 2005, Galucio et al. 2006, Poosh et al. 2012, Khosravian-Arab and Torres, 2013) tried to describe the properties of a discrete fractional filter in the time domain with causal functions. Then there are always transient effects at the initial time $t = 0$. There the filter has other properties than elsewhere because of the discontinuity of the input signal (and the resulting high frequencies).
When using the transfer function the realms of lower and of higher frequencies need special attention. For the lower frequencies the transfer function does not go to zero when using a finite number of points (in practice this is always the case). This we call the minimum frequency effect. When the input signal is causal the transfer function has no minimum frequency effect (this is important when for example treating differential equations). For the higher frequencies there is always a maximum depending on the sample frequency (Shannon frequency). But for the orthogonal derivative there is an extra maximum that is lower then the Shannon frequency. For the GL filter there is no such maximum.
A first conclusion of this section is that when using a discrete fractional filter the properties of the filter should be examined using the transfer function. A second conclusion is that a discrete fractional differentiating orthogonal filter has always a certain bandwidth.
In general a fractional differentiating filter can be built from two serial filters. The first one has a transfer function with modulus $|\omega|^v$. The second filter is a so-called low-pass filter. This low-pass filter determines the highest frequency for differentiation of the filter. In the case of the Jacobi filter (continuous case) one has to keep in mind that the transfer function has side lobes. In the discrete case these side lobes are not important because there is a maximum sample frequency. The Jacobi filter can be used if the maximum signal frequency is far beyond the side lobes. These sides lobes go to zero for high frequencies. Hence this is a very good fractional differentiating filter.
As another example for the low-pass filter we can choose the Butterworth filter. For this filter the transfer function is

$$\text{(11.20)} \quad H_n(\omega) = (-i\omega)^v \left( \frac{1}{1 + (\omega/\omega_0)^{2n}} \right)$$

which gives a very good transfer function for fractional differentiation. There are no side lobes. In Figure 7 this is demonstrated. There is no simple discrete analogue of this filter.
Diekema and Koornwinder (2012) make some remarks about the practical applications of the filters obtained from orthogonal polynomials and the Butterworth filter. In the analog integer case the Butterworth filter can be much easier constructed physically. In the analog fractional case a filter should be built with the transfer function $(-i\omega)^v$. This gives a difficult problem. Solutions are always approximations.

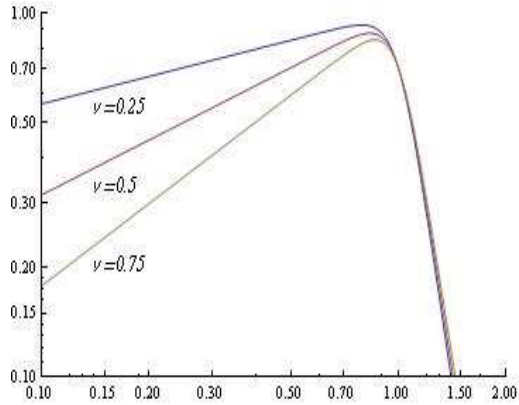

Figure 7. Modulus of the transfer function for the fractional Butterworth derivative with $n = 7$.

In the discrete case the Butterworth transfer function can be transformed into a filter equation. See Hamming (1989). For more details about the Butterworth filter see Oppenheim and Schafer (1975), §5.2.1, (both analog as discrete), Ziemer et al. (1998) (both analog as discrete), Johnson (1976), §3.2, (analog) and Hamming (1989), §12.6, (digital).

**Acknowledgments**

I thank Prof. Dr. T.H. Koornwinder for his help and time offered during my writing of this paper. Especially his hints about the Hahn polynomials were very helpful. Without his very stimulating enthusiasm and everlasting patience this work could not be done.

**Appendix A**. **The inverse Fourier transform of the function** $F(\omega) = (i\omega + 0)^{b-1}M(a,c;i\omega)$.

In this section we prove

(A.1) $\quad \mathcal{F}^{-1}\left[(i\omega + 0)^{b-1}M(a,c;i\omega)\right](y) =$

$$= \begin{cases} 0 & 1 < y \\ \sqrt{2\pi}\,\dfrac{\Gamma(c)}{\Gamma(a)\Gamma(1-a-b+c)} y^{a-b}(1-y)^{c-a-b} F\!\left(\begin{array}{c}1-b,c-b\\c+1-a-b\end{array};1-y\right) & 0 < y < 1 \\ \dfrac{\sqrt{2\pi}}{\Gamma(1-b)}\,\dfrac{1}{(1-y)^b} F\!\left(\begin{array}{c}b,c-a\\c\end{array};\dfrac{1}{1-y}\right) & y < 0 \end{cases}$$

with conditions for the parameters: $a,b,c \in \mathbb{R}$ and $0 < b < \min(a, c-a)$.

The inverse Fourier transform is given by

(A.2) $\quad I = \mathcal{F}^{-1}\left[(i\omega + 0)^{b-1}M(a,c;i\omega)\right](y) = \dfrac{1}{\sqrt{2\pi}}\int_{-\infty}^{\infty}(i\omega + 0)^{b-1}M(a,c;i\omega)e^{-i\omega y}d\omega =$

$\quad = \dfrac{e^{i\pi(b-1)/2}}{\sqrt{2\pi}}\int_{0}^{\infty}\omega^{b-1}M(a,c;i\omega)e^{-i\omega y}d\omega + \dfrac{e^{-i\pi(b-1)/2}}{\sqrt{2\pi}}\int_{0}^{\infty}\omega^{b-1}M(a,c;-i\omega)e^{i\omega y}d\omega$

If $a,c \in \mathbb{R}$ then $|M(a;c;ix)| \leq C|x|^{\max(-a,a-c)}$ [Erdélyi (1953) formula 6.13(2)]. The conditions for the parameters that the above integrals converge absolutely are $0 < b < \min(a, c-a)$, so the inverse Fourier transform exists. We can rewrite the integral in (A.2) as

(A.3) $\quad I = \dfrac{e^{i\pi(b-1)/2}}{\sqrt{2\pi}}\lim_{\varepsilon\downarrow 0}\int_{0}^{\infty}\omega^{b-1}M(a,c;i\omega)e^{-\omega(\varepsilon+iy)}d\omega + \dfrac{e^{-i\pi(b-1)/2}}{\sqrt{2\pi}}\lim_{\varepsilon\downarrow 0}\int_{0}^{\infty}\omega^{b-1}M(a,c;-i\omega)e^{-\omega(\varepsilon-iy)}d\omega$

The integrals are known [Erdélyi (1953) formula 6.10(5)]. There follows

(A.4) $$I = \frac{\Gamma(b)}{\sqrt{2\pi}} \lim_{\varepsilon \downarrow 0} \frac{e^{i\pi(b-1)/2}}{(iy+\varepsilon)^b} F\left(\begin{matrix}a,b\\c\end{matrix}; \frac{i}{\varepsilon+iy}\right) + \frac{\Gamma(b)}{\sqrt{2\pi}} \lim_{\varepsilon \downarrow 0} \frac{e^{-i\pi(b-1)/2}}{(\varepsilon-iy)^b} F\left(\begin{matrix}a,b\\c\end{matrix}; \frac{-i}{\varepsilon-iy}\right)$$

Taking the limit gives

(A.5) $$I = -i\frac{\Gamma(b)}{\sqrt{2\pi}} \frac{1}{(y-i0)^b} F\left(\begin{matrix}a,b\\c\end{matrix}; \frac{1}{y-i0}\right) + i\frac{\Gamma(b)}{\sqrt{2\pi}} \frac{1}{(y+i0)^b} F\left(\begin{matrix}a,b\\c\end{matrix}; \frac{1}{y+i0}\right)$$

Because there is a cut we distinguish three possibilities.

I. For $y > 1$ there remains

(A.6) $$I = i\frac{\Gamma(b)}{\sqrt{2\pi}} \frac{1}{y^b} F\left(\begin{matrix}a,b\\c\end{matrix}; \frac{1}{y}\right) - i\frac{\Gamma(b)}{\sqrt{2\pi}} \frac{1}{y^b} F\left(\begin{matrix}a,b\\c\end{matrix}; \frac{1}{y}\right) = 0$$

II. For $0 < y < 1$ there remains

(A.7) $$I = i\frac{\Gamma(b)}{\sqrt{2\pi}} \frac{1}{(y+i0)^b} F\left(\begin{matrix}a,b\\c\end{matrix}; \frac{1}{y+i0}\right) - i\frac{\Gamma(b)}{\sqrt{2\pi}} \frac{1}{(y-i0)^b} F\left(\begin{matrix}a,b\\c\end{matrix}; \frac{1}{y-i0}\right)$$

We use formula 2.10(4) of Erdélyi (1953). These formula gives

(A.8) $$F\left(\begin{matrix}a,b\\c\end{matrix}; z\right) = \frac{\Gamma(c)\Gamma(c-a-b)}{\Gamma(c-a)\Gamma(c-b)} \frac{1}{z^a} F\left(\begin{matrix}a, a+1-c\\a+b+1-c\end{matrix}; 1-\frac{1}{z}\right) +$$
$$+ \frac{\Gamma(c)\Gamma(a+b-c)}{\Gamma(a)\Gamma(b)} \frac{(1-z)^{c-a-b}}{z^{c-a}} F\left(\begin{matrix}c-a, 1-a\\c+1-a-b\end{matrix}; 1-\frac{1}{z}\right)$$

Application gives

(A.9) $$I = i\frac{\Gamma(b)}{\sqrt{2\pi}} \frac{1}{y^b} \frac{\Gamma(c)\Gamma(a+b-c)}{\Gamma(a)\Gamma(b)} y^{c-a} \left(1 - \frac{1}{y+i0}\right)^{c-a-b} F\left(\begin{matrix}c-a, 1-a\\c+1-a-b\end{matrix}; 1-y\right) -$$
$$- i\frac{\Gamma(b)}{\sqrt{2\pi}} \frac{1}{y^b} \frac{\Gamma(c)\Gamma(a+b-c)}{\Gamma(a)\Gamma(b)} y^{c-a} \left(1 - \frac{1}{y-i0}\right)^{c-a-b} F\left(\begin{matrix}c-a, 1-a\\c+1-a-b\end{matrix}; 1-y\right)$$

After a little computation there follows

(A.10) $$I = \sqrt{2\pi} \frac{\Gamma(c)}{\Gamma(a)\Gamma(1-a-b+c)} (1-y)^{c-a-b} F\left(\begin{matrix}c-a, 1-a\\c+1-a-b\end{matrix}; 1-y\right)$$

Application of Olver et al. (2010) formula 15.8. gives at last

(A.11) $$I = \sqrt{2\pi} \frac{\Gamma(c)}{\Gamma(a)\Gamma(1-a-b+c)} y^{a-b}(1-y)^{c-a-b} F\left(\begin{matrix}1-b, c-b\\c+1-a-b\end{matrix}; 1-y\right)$$

III. For $y < 0$ there remains

(A.12) $$I = i\frac{\Gamma(b)}{\sqrt{2\pi}} \frac{1}{(-y)^b} \left[e^{-i\pi b} - e^{i\pi b}\right] F\left(\begin{matrix}a,b\\c\end{matrix}; \frac{1}{y}\right) = \frac{\sqrt{2\pi}}{\Gamma(1-b)} \frac{1}{(-y)^b} F\left(\begin{matrix}a,b\\c\end{matrix}; \frac{1}{y}\right)$$

Application of Olver et al. (2010) formula 15.8.1 gives

(A.13) $$I = \frac{\sqrt{2\pi}}{\Gamma(1-b)} \frac{1}{(1-y)^b} F\left(\begin{matrix}b, c-a\\c\end{matrix}; \frac{1}{1-y}\right)$$

Gathering the results proves (A.1).□